\documentclass[12pt]{article}

\usepackage{latexsym}
   \usepackage[all]{xy}
   \usepackage{amsfonts}
   \usepackage{amsthm}
   \usepackage{amsmath}
   \usepackage{amssymb}
\numberwithin{equation}{section}
   \usepackage[all]{xy}
\usepackage{mathptmx}

\def\<{\langle}
\def\>{\rangle}

\def\tens{\mathop{\otimes}}
\def\la{{\triangleright}}\def\ra{{\triangleleft}}

\def\extd{{\rm d}}

\def\id{{\rm id}}

   \newcommand{\cC}{{\mathfrak C}}
     \def \eC{{\epsilon_\cC}}
       \def\ut{{\otimes}}
         \def\DC{{\Delta_\cC}}
           \def\ot{{\otimes}}
             \def\cD{{\mathfrak D}}
    \def \eD{{\epsilon_\cD}}
   \def\sw#1{{\sb{(#1)}}}
    \def\sco#1{{\sb{[#1]}}}
     \def\M{\mathcal{M}}
       \def\DD{{\Delta_\cD}}
       \def\db{\mathsf{DiffBim}}
       \def\fdb{\mathsf{FlatDiffBim}}

\newcommand{\LEM}{\mathsf{LEM}(\mathsf{Bim})}
\def\tto{\longrightarrow}
\def\Z{\mathbb{Z}}
\def\deg{\mathrm{deg}\,}
	
\newtheorem{lemma}{Lemma}[section]
\newtheorem{propos}[lemma]{Proposition}
\newtheorem{theorem}[lemma]{Theorem}
\newtheorem{cor}[lemma]{Corollary}

\theoremstyle{definition}
\newtheorem{defin}[lemma]{Definition}
\newtheorem{example}[lemma]{Example}

\theoremstyle{remark}
\newtheorem{remark}[lemma]{Remark}

\begin{document}
\title{The Serre spectral sequence of a noncommutative fibration
for de~Rham cohomology}
\author{E.J.\ Beggs \& Tomasz Brzezi\'nski \\  \\ Department of 
Mathematics\\ University of Wales Swansea \\
Swansea SA2 8PP, U.K.}
\date{August 2005}

\maketitle

\abstract{For differential calculi on noncommutative algebras, we construct
a twisted de Rham cohomology using flat connections on modules. This 
has properties similar, in some respects, to sheaf cohomology on 
topological spaces. We also discuss generalised mapping properties of 
these theories, and relations of these properties to corings.
Using this, we
give conditions for the Serre spectral sequence to hold for a noncommutative
fibration. This might be better read as giving the definition of a fibration in
noncommutative differential geometry. We also study the 
multiplicative structure of
such spectral sequences. Finally we show that some noncommutative 
homogeneous spaces satisfy the conditions to be such a fibration, and 
in the process
clarify the differential structure on these homogeneous spaces. We 
also give two explicit examples of differential fibrations: these are 
built on the quantum Hopf fibration with two different differential 
structures.}

\section{Introduction}
This paper has three basic purposes:

\smallskip \noindent
1) Developing a cohomology theory for modules with flat 
connections
over noncommutative algebras, and showing that it has 
some 
properties in common with sheaf theory. 

 \noindent
2)  Extending the Serre spectral sequence of a fibration in classical algebraic
topology to the noncommutative domain.

\noindent
3) Examining the differential structure of quantum homogeneous 
spaces, and showing that
many of them are `fibrations' in a noncommutative sense.

\smallskip

In \cite{ConnesNCDG} methods of studying algebras by means of their differential
calculi were introduced. We will apply Connes' differential methods 
to fibrations
in algebraic topology.

In usual topology, sheaf cohomology and other methods
allow cohomology with `twisted' coefficients, i.e.\ coefficients
which vary from point to point in the space.
  In the absence (so far at least)
of a full sheaf cohomology construction in noncommutative geometry,
we construct de Rham cohomology with twisted coefficients
for algebras with differential structure. The allowed coefficients are
modules with flat connection. Though
there is a considerable similarity between de Rham cohomology with twisted
coefficients in the noncommutative world
  and sheaf cohomology in the commutative world, it is
quite possible that yet more general constructions, or constructions 
with additional properties,
corresponding to sheaf theory exist in the noncommutative world.
In the spirit of some developments in operator algebra (for example, see
\cite{EvaKaw:book}), we show
that bimodules can be used to replace algebra maps in
constructing `pull backs' of the coefficient modules. In the special case of
semi-free differential graded algebras this construction is shown to
have an interesting interpretation in terms of {\em corings}.

In commutative algebraic topology, one of the most useful
applications of twisted coefficients is to fibrations. For a locally trivial
fibration, the Serre spectral sequence \cite{spseq} starts with the 
cohomology of
the base space with coefficients in the cohomology of the fibre
(in general a twisted bundle), and converges to the cohomology of the
total space. In producing a noncommutative analogue of this result,
we not only have to find a proof which does not require local triviality,
but also have to decide what a `locally trivial' fibration should be
in noncommutative differential geometry. Realistically we should define a
  fibration by the conditions which are required by the Serre spectral sequence.
  The seeming correspondence between sheaf theory and the cohomology we
  are considering leads us to suspect that yet more general `Leray type'
  spectral sequences exist.

We then discuss products in the Serre
spectral sequence, which requires another condition to be imposed on the
fibration. The product structure is not only important in its own right, but
can frequently help in simplfying the calculation of  spectral sequences.

Finally, we
  study fibrations given in terms of a coaction of a Hopf algebra on an algebra.
  As a non-trivial example of such a differential fibration we 
consider the quantum Hopf fibration $\iota: 
\mathcal{A}(S_q^2)\hookrightarrow \mathcal{A}(SL_q(2))$ with the 
three-dimensional differential calculus  on $\mathcal{A}(SL_q(2))$. 
As a further non-trivial class of examples of the fibrations 
discussed here, we look at
the noncommutative homogeneous space construction
 with bicovariant 
differential calculi. This takes the
classical construction of a group quotiented by a subgroup, and replaces it by
two Hopf algebras with a surjective Hopf algebra map $\pi:X\to H$.
We begin with such a $\pi$ which is differentiable with respect to
bicovariant differential structures on $X$ and $H$ \cite{worondiff}. Note that
the bicovariant condition corresponds to the differentiability of the
coproduct map, and it is reasonable to expect that this is
the analogue of classical Lie groups. As in the classical case, some 
of the definitions can be given in terms of the Hopf-Lie algebras and 
their induced vector fields
 \cite{braexp}. 
   The first stage is to 
identify the differential calculus for the homogeneous space 
$B=X^{{\rm co}H}$ (see Theorem
\ref{bchsinbjic}) in a form suitable for calculation.
Then it is shown that
the inclusion map $B \to X$ is a fibration as defined earlier (see 
Theorem \ref{ikvyu}).
For the development of
noncommutative homogenous spaces
  the reader should refer to \cite{HecSch,SMspin}. Again the quantum 
Hopf fibration
  $\iota: \mathcal{A}(S_q^2)\hookrightarrow \mathcal{A}(SL_q(2))$ is 
an example of this situation and we explicitly prove that it is a 
differentiable fibration for one of two standard four-dimensional 
bicovariant calculi on $\mathcal{A}(SL_q(2))$.

All algebras are unital over a field $k$. The unadorned
tensor product between vector spaces is over $k$. A Hopf algebra is 
always assumed
to have a bijective antipode (this is not the most general situation 
algebraically,
but the most natural from the point of view of non-commutative
geometry).

\section{Flat connections and cohomology with twisted coefficients}
The classical Serre spectral sequence uses cohomology with a 
non-trivial coefficient
bundle. In this section we discuss flat connections on modules in 
noncommutative geometry, and how this can be used to define de Rham 
cohomology with non-trivial
coefficient modules.

By a {\em differential calculus} on a noncommutative algebra $A$ we mean
a differential graded algebra $({\rm d},\Omega^* A)$ such that $\Omega^0A = A$.
The product in $\Omega^* A$ (for $*\geq 1$) is denoted by the wedge 
$\wedge$ (although $\Omega^* A$ is not graded anticommutative in 
general). The density condition says that
$\Omega^{n+1} A\subset A.\extd\Omega^n A$, but we will not require 
this till later.

The cohomology of $({\rm d},\Omega^* A)$ is denoted by $H_{dR}^*(A)$ and
referred to as a {\em de Rham cohomology} of $A$. Recall that a {\em 
connection}
in a left $A$-module $E$ is a map $\nabla: E\to\Omega^1 A \tens_A E$ satisfying
the Leibniz rule, for all $a\in A$, $e\in E$, $\nabla (a.e) = \extd 
a\otimes e + a\nabla e$.

\subsection{The construction of the cohomology}

\begin{defin}
Given an algebra $A$ with differential calculus $({\rm d},\Omega^* A)$,
we define the category ${}_A\mathcal{E}$ to consist of
left $A$-modules $E$ with connection $\nabla:E\to\Omega^1 A \tens_A E$.
A morphism $\phi:(E,\nabla)\to (F,\nabla)$  in the category is a left 
$A$-module map
$\phi:E\to F$ which preserves the covariant derivative, i.e.
$\nabla\circ\phi=(\id\tens\phi)\circ\nabla:E\to \Omega^1 A \tens_A F$.
\end{defin}

\begin{defin}
Given $(E,\nabla)\in {}_A\mathcal{E}$, define 
$$\nabla^{[n]}:\Omega^n A
\tens_A E\to \Omega^{n+1} A\tens_A E\, , \qquad \omega\tens e
\mapsto 
\extd\omega\tens e+(-1)^n\omega\wedge \nabla e. 
$$
Then the {\em curvature}
is defined as $R=\nabla^{[1]}\nabla:E\to \Omega^2 A\tens_A E$, and is a left
$A$-module map.

  The covariant derivative is called {\em flat}
if the curvature is zero. We write ${}_A\mathcal{F}$ for the full
subcategory of ${}_A\mathcal{E}$
consisting of left $A$-modules with flat connections.
\end{defin}

\begin{propos}
For all $n\ge 0$, $\nabla^{[n+1]}\circ \nabla^{[n]}=\id\wedge R:\Omega^n A
\tens_A E\to \Omega^{n+2} A\tens_A E$.
\end{propos}
\noindent
{\bf Proof:}\quad By explicit calculation,
\begin{eqnarray*}
\nabla^{[n+1]}( \nabla^{[n]}(\omega\tens e)) &=&
\nabla^{[n+1]}(\extd\omega\tens e+(-1)^n\omega\wedge \nabla e)\ .
\end{eqnarray*}
Put $\nabla e=\xi_i\tens e_i$ (summation implicit), and then
\begin{eqnarray*}
\nabla^{[n+1]}( \nabla^{[n]}(\omega\tens e)) &=&
\nabla^{[n+1]}(\extd\omega\tens e+(-1)^n\omega\wedge \xi_i\tens e_i) \cr
&=& (-1)^{n+1}\extd\omega\wedge \nabla e+
(-1)^n\extd\omega\wedge \xi_i\tens e_i
+\omega\wedge \extd\xi_i\tens e_i \cr
&& -\, \omega\wedge \xi_i\wedge \nabla e_i \cr
&=& \omega\wedge(\extd\xi_i\tens e_i
- \xi_i\wedge \nabla e_i)\,=\,\omega\wedge R(e)\ .\quad\square
\end{eqnarray*}

\begin{defin} \label{codef}
Given $(E,\nabla)\in {}_A\mathcal{F}$, define $H^*(A;E,\nabla)$ to be the
cohomology of the complex
\begin{eqnarray*}
E \stackrel{\nabla}{\longrightarrow}
\Omega^1 A\tens_A E  \stackrel{\nabla^{[1]}}{\longrightarrow}
\Omega^2 A\tens_A E  \stackrel{\nabla^{[2]}}{\longrightarrow} \dots\ .
\end{eqnarray*}
Note that $H^0(E,\nabla)=\Gamma E=\{e\in E:\nabla e=0\}$, the flat 
sections of $E$.
We will often write $H^*(A;E)$ where there is no danger of confusing
the covariant derivative.
\end{defin}

\begin{propos} \label{kfen}
Given $(E,\nabla)\in {}_A\mathcal{F}$,
the map $\wedge:\Omega^n A\tens(\Omega^r A\tens_A E)
\to \Omega^{n+r} A\tens_A E$ defined by
$\wedge(\xi\tens(\omega\tens e))=(\xi\wedge \omega)\tens e$
gives a graded left $H_{dR}(A)$-module structure on $H^*(A;E,\nabla)$.
\end{propos}
\noindent {\bf Proof:}\quad First calculate
\begin{eqnarray*}
\nabla^{[*]}(\xi\wedge(\omega\tens e)) &=& 
\nabla^{[*]}((\xi\wedge\omega)\tens e) \cr
&=& \extd (\xi\wedge\omega)\tens e +(-1)^{|\xi|+|\omega|}
\xi\wedge\omega\wedge\nabla e \cr
&=& \extd\xi\wedge (\omega\tens e)+(-1)^{|\xi|}\xi\wedge\nabla^{[*]}
  (\omega\tens e)\ .
\end{eqnarray*}
This equation has the required immediate consequences:

If $\extd\xi=0$ and $\nabla^{[*]}(\omega\tens e)=0$ then
$\nabla^{[*]}(\xi\wedge(\omega\tens e))=0$.

If $\nabla^{[*]}(\omega\tens e)=0$
then $\extd\xi\wedge (\omega\tens e)$ is in the image of $\nabla^{[*]}$.

If $\extd\xi=0$
then $\xi\wedge \nabla^{[*]}(\omega\tens e)$ is in the image of $\nabla^{[*]}$.
\quad$\square$

\subsection{Mapping properties of the cohomology}
In classical topology, maps on the cohomology can
be induced by maps which change coefficients over the same topological space.
Our analogue of this is the following:

\begin{theorem}
The cohomology $H^*$ in Definition~\ref{codef} is a functor from
${}_A\mathcal{F}$ to graded left $H_{dR}^*(A)$-modules,
where the module structure is given in Proposition~\ref{kfen}.
\end{theorem}
\noindent
{\bf Proof:}\quad Begin with a left $A$-module map
$\phi:E\to F$ which preserves the covariant derivative, i.e.
$\nabla\circ\phi=(\id\tens\phi)\circ\nabla:E\to \Omega^1 A \tens_A F$.
First show that the map $\id\tens\phi:
\Omega^*A\tens_A E\to \Omega^*A\tens_A F$ is a cochain map:
\begin{eqnarray*}
\nabla^{[*]}(\id\tens\phi)(\omega\tens e) &=&
\nabla^{[*]}(\omega\tens\phi(e)) \cr
&=& \extd\omega\tens\phi(e) +(-1)^{|\omega|}\,\omega\wedge\nabla\phi(e) \cr
&=& \extd\omega\tens\phi(e) +(-1)^{|\omega|}\,\omega\wedge(\id\tens\phi)
\nabla e \cr
&=& (\id\tens \phi)\nabla^{[*]}(\omega\tens e)\ .
\end{eqnarray*}
The functorial property is simply $(\id\tens\phi)\circ(\id\tens\psi)=
\id\tens(\phi\circ\psi)$, and the left module property is just
$\xi\wedge(\omega\tens\phi(e))=(\xi\wedge\omega)\tens\phi(e)$. \quad$\square$

\medskip In classical topology, continuous functions between topological spaces
also induce maps on the cohomology. One part of this is the pull back 
construction
for coefficients. Given the reversal of arrows which often occurs in 
considering algebras rather than spaces, this becomes a `push 
forward' construction in non-commutative geometry.

This would be an appropriate time to remind the reader that for 
algebras $A$ and $B$
with differentiable structure, an algebra map $\theta:A\to B$ is 
called differentiable
if it extends to a map $\theta_*:\Omega^*A\to\Omega^* B$ of 
differential graded algebras.

\begin{lemma} \label{bjhksa}
Given $(E,\nabla)\in {}_A\mathcal{E}$ and a differentiable algebra map
$\theta:A\to B$,
define 
$$\hat\nabla:B\tens_A E\to \Omega^1 B\tens_B 
B\tens_A E
=\Omega^1 B\tens_A E \, , \quad
  \hat\nabla(b\tens e)=b.(\theta_*\tens\id)(\nabla e)+\extd b\tens e.
  $$
  Then $\theta_*(E,\nabla)=(B\tens_A E,\hat\nabla)\in 
{}_B\mathcal{E}$, with right action
of $A$ on $B$ given by
$b\ra a=b\,\theta(a)$.
\end{lemma}
\noindent
{\bf Proof:} To check that $\hat\nabla$ is well defined,
we must show that, for all $a\in A$, $b\in B$ and $e\in E$,
$\hat\nabla(b\,\theta(a)\tens e)=\hat\nabla(b\tens a\,e)$:
\begin{eqnarray*}
\hat\nabla(b\,\theta(a)\tens e) &=& b\,\theta(a).
(\theta_*\tens\id)(\nabla e)+\extd (b\,\theta(a))\tens e \cr
&=& b.(\theta_*\tens\id)(a\,\nabla e)+\extd b.\theta(a)\tens e+
b.\extd\theta(a)\tens e \cr
&=& b.(\theta_*\tens\id)(a\,\nabla e+\extd a\tens e)+\extd b\tens a\,e \cr
&=& b.(\theta_*\tens\id)\nabla(a\, e)+\extd b\tens a\,e 
= \hat\nabla(b\tens a\,e)\ .
\end{eqnarray*}
That $\hat\nabla$ satisfies the Leibniz rule follows immediately from the 
definition (and the Leibniz rule for $\extd$). $\quad\square$

\begin{propos} \label{vcbsdhiu}
If $\theta:A\to B$ is a differentiable algebra map and
$(E,\nabla)\in {}_A\mathcal{F}$, then $\theta_*(E,\nabla)\in {}_B\mathcal{F}$.
\end{propos}
\noindent {\bf Proof:}\quad
Following the notation of Lemma~\ref{bjhksa} and setting $\nabla 
e=\xi_i\tens e_i$
(summation implied),
\begin{eqnarray*}
\hat\nabla^{[1]}\hat\nabla(b\tens e) &=&
\hat\nabla^{[1]}(b.(\theta_*\tens\id)(\nabla e)+\extd b\tens e) \cr
&=& \hat\nabla^{[1]}(b.\xi_i\tens e_i+\extd b\tens e) \cr
&=& \extd(b.\xi_i)\tens e_i+\extd\extd b\tens e-b.\xi_i\wedge\nabla e_i
-\extd b\wedge\nabla e \cr
&=& \extd b\wedge\xi_i\tens e_i+b.\extd \xi_i\tens e_i-b.\xi_i\wedge\nabla e_i
-\extd b\wedge\nabla e \cr
&=& b.(\extd \xi_i\tens e_i-\xi_i\wedge\nabla e_i) \cr
&=& b.\nabla^{[1]}\nabla e\,=\,0\ .\quad\square
\end{eqnarray*}

\begin{theorem} \label{hjhjhj}
For $\theta:A\to B$ a differentiable algebra map,
there is a functor $\theta_*:{}_A\mathcal{E}\to {}_B\mathcal{E}$
which is defined on objects as  in Lemma~\ref{bjhksa}, and where
a morphism $\phi:E\to F$ is sent to the morphism $\id\tens\phi:B\tens_A E
\to B \tens_A F$. Further this functor restricts to
a functor from ${}_A\mathcal{F}$ to ${}_B\mathcal{F}$.
\end{theorem}
\noindent
{\bf Proof:}\quad First, given a morphism $\phi:E\to F$ in ${}_A\mathcal{E}$,
we need to show that $\id\tens\phi:B\tens_A E
\to B \tens_A F$ is a morphism in ${}_B\mathcal{E}$. Using the definition
of $\hat\nabla$ in Lemma~\ref{bjhksa},
\begin{eqnarray*}
\hat\nabla(b\tens\phi(e)) &=& b.(\theta_*\tens\id)\nabla 
\phi(e)+\extd b\tens \phi(e)\ ,
\end{eqnarray*}
and as $\phi$ is a morphism in ${}_A\mathcal{E}$,
\begin{eqnarray*}
\hat\nabla(b\tens\phi(e)) &=& b.(\theta_*\tens\id)(\id\tens\phi)
\nabla e+\extd b\tens \phi(e)\cr
&=& (\id\tens\phi)(b.(\theta_*\tens\id)\nabla e+\extd b\tens e)\cr
&=& (\id\tens\phi)\hat\nabla(b\tens e)\ .
\end{eqnarray*}
The composition rule is just 
$(\id\tens\phi)\circ(\id\tens\psi)=\id\tens\phi\circ\psi$.
The restriction to flat connections is shown in 
Proposition~\ref{vcbsdhiu}.\quad$\square$

\subsection{Generalised mapping properties}
The mapping constructions can be generalised to bimodules rather
than algebra maps, using the `braiding' introduced by Madore \cite{Madore}:

\begin{defin} \label{tenten}
A $B$-$A$ bimodule $M\in{}_B \mathcal{M}_A$ with additional
structures

(a)\quad a left $B$-connection $\check \nabla:M\to \Omega^1 B\tens_B M$,

(b)\quad a $B$-$A$ bimodule map $\check\sigma:M\tens_A \Omega^1 A \to
\Omega^1 B\tens_B M$,

\noindent is called a {\em differentiable bimodule} if it satisfies the
condition $\check \nabla(m.a)=\check \nabla(m).a+\check\sigma(m\tens\extd a)$
for all $m\in M$ and $a\in A$.
\end{defin}

\begin{example}
If $\theta:A\to B$ is a differentiable algebra map,
take the bimodule $B\in {}_B \mathcal{M}_A$, with the usual
left $B$-action, and right $A$-action given by $b\ra a=b\,\theta(a)$.
Also define $\check\nabla:B\to \Omega^1 B \tens_B B=\Omega^1 B$
by $\check\nabla b=\extd b$ and $\check\sigma:B\tens_A \Omega^1 A \to
\Omega^1 B\tens_B B=\Omega^1 B$ by $\check\sigma(b\tens\xi)
=b.\theta_*(\xi)$. Now we check the condition
\begin{eqnarray*}
\check \nabla(b\ra a) &=& \check \nabla(b\,\theta(a))\,=\,\extd 
(b\,\theta(a)) \cr
&=& \extd b\,\theta(a) + b\,\theta_*(\extd a)\,=\,
\check \nabla(b).a+\check\sigma(b\tens\extd a)\ .
\end{eqnarray*}
Hence $B$ is a differentiable bimodule.
\end{example}

\begin{propos} \label{prop.functor} Suppose that
  $(M,\check \nabla,\check\sigma)$ is a differentiable $B$-$A$ bimodule.
   Then the following
defines a functor $(M,\check \nabla,\check\sigma)_*:
{}_A\mathcal{E}\to {}_B\mathcal{E}$:

On objects $(E,\nabla)\in {}_A\mathcal{E}$, define
$(M,\check \nabla,\check\sigma)_*(E,\nabla)=(M\tens_A E,\hat\nabla)$,
where
\begin{eqnarray*}
\hat\nabla(m\tens e)\,=\,\check\nabla m\tens e+
(\check\sigma\tens\id)(m\tens\nabla e)
\end{eqnarray*}

On morphisms $\phi:E\to F$, define $(M,\check \nabla,\check\sigma)_*\phi=
\id\tens\phi:M\tens_A E\to M\tens_A F$.
\end{propos}
\noindent
{\bf Proof:}\quad First we need to check that $\hat\nabla$ is a well defined
function on $M\tens_A E$.
\begin{eqnarray*}
\hat\nabla(m\tens a\, e) &=& \check\nabla m\tens a\, e+
(\check\sigma\tens\id)(m\tens\nabla (a\,e)) \cr
&=& (\check\nabla m).a\tens a e+
(\check\sigma\tens\id)(m\tens a\,\nabla e) 
+(\check\sigma\tens\id)(m\tens\extd a\tens e) \ .
\end{eqnarray*}
By using the differentiable bimodule condition this becomes
\begin{eqnarray*}
\hat\nabla(m\tens a\, e)
&=& \check\nabla (m.a)\tens a e+
(\check\sigma\tens\id)(m.a\tens \nabla e)  \,=\,\hat\nabla(m.a\tens  e) \ .
\end{eqnarray*}
To check that $\hat\nabla$ is a left $B$-covariant derivative,
as $\check\sigma$ is a left $B$-module map,
\begin{eqnarray*}
\hat\nabla(b.m\tens e) &=& \check\nabla (b.m)\tens e+
(\check\sigma\tens\id)(b.m\tens\nabla e) \cr
&=& b.\check\nabla (m)\tens e+\extd b\tens m\tens e +
b.(\check\sigma\tens\id)(m\tens\nabla e) \cr
&=& b.\hat\nabla(m\tens e)+\extd b\tens m\tens e\ .
\end{eqnarray*}
Next we check the morphism condition,
\begin{eqnarray*}
\hat\nabla(m\tens \phi(e)) &=&
\check\nabla m\tens \phi(e)+
(\check\sigma\tens\id)(m\tens\nabla \phi(e)) \cr &=&
\check\nabla m\tens \phi(e)+
(\check\sigma\tens\phi)(m\tens\nabla e) \cr
&=&(\id\tens\id\tens\phi)\hat\nabla(m\tens e) \ .\quad\square
\end{eqnarray*}

\begin{defin} \label{bhcsial}
The differentiable $B$-$A$ bimodule  $(M,\check \nabla,\check\sigma)$ is
said to be  {\em flat} if there is a  $B$-$A$ bimodule map 
$\check\sigma:M\tens_A \Omega^2 A \to
\Omega^2 B\tens_B M$ so that the following conditions are satisfied:

(a)\quad as a left $B$-connection on $M$, $\check\nabla$ is flat;

(b)\quad $(\id\wedge\check\sigma)(\check\sigma\tens\id) =
\check\sigma(\id\tens\wedge):M\tens_A\Omega^1 A\tens_A\Omega^1 A\to
\Omega^2 B\tens_B M$.

\end{defin}

\begin{lemma} \label{lkvbihdz}
If the differentiable 
$B$-$A$ bimodule  $(M,\check \nabla,\check\sigma)$ is {\em 
flat},
then the following map 
vanishes:
\begin{eqnarray*}
[(\extd\tens\id)
-(\id\wedge\check\nabla)]
\check\sigma -(\id\wedge\check\sigma)(\check\nabla\tens\id)
-\check\sigma (\id\tens\extd):M\tens_A\Omega^1 A\to \Omega^2 B\tens_B 
M
\end{eqnarray*}
\end{lemma}
\noindent
{\bf Proof:}\quad First note that the displayed formula is well 
defined,
as for all $m\in M$, $\eta\in\Omega^1 B$ and $b\in 
B$;
\begin{eqnarray*}
[(\extd\tens\id)-(\id\wedge\check\nabla)](\eta\,b\tens 
m) \,=\,
[(\extd\tens\id)-(\id\wedge\check\nabla)](\eta\tens b\,m) \ 
.
\end{eqnarray*}
Next we check that the displayed formula defines a 
right $A$-module map. 
For all $m\in M$, $a\in A$ and $\xi\in\Omega^1 
A$,
\begin{eqnarray*}
[(\extd\tens\id)-(\id\wedge\check\nabla)]\check\sigma(m\tens\xi\,a) 
&=&
[(\extd\tens\id)-(\id\wedge\check\nabla)](\check\sigma(m\tens\xi)\,a) 
\cr
&=&
[(\extd\tens\id)-(\id\wedge\check\nabla)](\check\sigma(m\tens\xi))\,a 
\cr
&& -\,  (\id\wedge 
\check\sigma)(\check\sigma(m\tens\xi)\tens\extd a)\ 
,\cr
[(\id\wedge\check\sigma)(\check\nabla\tens\id)
+\check\sigma (\id\tens\extd)](m\tens\xi\,a) &=& 
[(\id\wedge\check\sigma)(\check\nabla\tens\id)
+\check\sigma (\id\tens\extd)](m\tens\xi)\,a \cr
&&-\, \check\sigma 
(m \tens \xi \wedge\extd a)\ ,
\end{eqnarray*}
and combining these 
and using Definition \ref{bhcsial} gives right $A$-linearity.

Finally to prove the vanishing of the displayed formula, we now 
only have to apply it to
elements of the form $m\tens \extd a$, and 
use the differentiable bimodule condition on 
$\check\sigma$.
\begin{eqnarray*}
[(\extd\tens\id)-(\id\wedge\check\nabla)]\check\sigma(m\tens\extd 
a) 
&=&
[(\extd\tens\id)-(\id\wedge\check\nabla)](\check\nabla(m.a)-(\check\nabla 
m).a)\cr
&=& 
[(\extd\tens\id)-(\id\wedge\check\nabla)]\check\nabla(m.a) \cr
&&-\, 
[(\extd\tens\id)-(\id\wedge\check\nabla)]\check\nabla(m).a \cr
&& +\, 
(\id\wedge\check\sigma)(\check\nabla m\tens \extd a)\ 
,\cr
[(\id\wedge\check\sigma)(\check\nabla\tens\id)
+\check\sigma (\id\tens\extd)](m\tens\extd a) 
&=&
(\id\wedge\check\sigma)(\check\nabla m\tens \extd a)\ 
.
\end{eqnarray*}
This means that the displayed formula applied to 
$m\tens \extd a$ gives
$R(m.a)-R(m).a$, where $R$ is the curvature of 
the left $B$-connection on $M$,
and this vanishes by Definition 
\ref{bhcsial}. \quad$\square$

\begin{propos} \label{fourteen}
If the differentiable $B$-$A$ bimodule  $(M,\check 
\nabla,\check\sigma)$ is flat,
then the functor $(M,\check \nabla,\check\sigma)_*:
  {}_A\mathcal{E}\to  {}_B\mathcal{E}$ restricts to a functor
from $ {}_A\mathcal{F}$ to $ {}_B\mathcal{F}$.
\end{propos}
\noindent {\bf Proof:}\quad We need to show that the following vanishes:
\begin{eqnarray*}
\hat\nabla^{[1]}\hat\nabla(m\tens e) &=& 
\hat\nabla^{[1]}(\check\nabla m\tens e+
(\check\sigma\tens\id)(m\tens\nabla e)) \cr
&=& (\extd\tens\id\tens\id)(\check\nabla m\tens e+
(\check\sigma\tens\id)(m\tens\nabla e)) \cr
&&-\,(\id\wedge\hat\nabla)(\check\nabla m\tens e+
(\check\sigma\tens\id)(m\tens\nabla e)) \cr
&=& (\extd\tens\id\tens\id)(\check\nabla m\tens e+
(\check\sigma\tens\id)(m\tens\nabla e)) \cr
&&-\,(\id\wedge\check\nabla\tens\id)(\check\nabla m\tens e+
(\check\sigma\tens\id)(m\tens\nabla e)) \cr
&&-\,(\id\wedge\check\sigma\tens\id)
(\id\tens\id\tens\nabla)(\check\nabla m\tens e+
(\check\sigma\tens\id)(m\tens\nabla e)) \cr
&=& (\extd\tens\id\tens\id)(\check\nabla m\tens e)+ (\extd\tens\id\tens\id)
(\check\sigma\tens\id)(m\tens\nabla e) \cr
&&-\,(\id\wedge\check\nabla\tens\id)(\check\nabla m\tens e)
-(\id\wedge\check\nabla\tens\id)
(\check\sigma\tens\id)(m\tens\nabla e) \cr
&&-\,(\id\wedge\check\sigma\tens\id)
(\id\tens\id\tens\nabla)(\check\nabla m\tens e+
(\check\sigma\tens\id)(m\tens\nabla e))\ .
\end{eqnarray*}
As the left $B$ covariant derivative $\check\nabla$ on $M$ is flat, 
the first and third terms
cancel, giving
\begin{eqnarray*}
\hat\nabla^{[1]}\hat\nabla(m\tens e)
&=&  (\extd\tens\id\tens\id)
(\check\sigma\tens\id)(m\tens\nabla e)
-(\id\wedge\check\nabla\tens\id)
(\check\sigma\tens\id)(m\tens\nabla e) \cr
&&-\,(\id\wedge\check\sigma\tens\id)
(\id\tens\id\tens\nabla)(\check\nabla m\tens e)\cr
&&-\,(\id\wedge\check\sigma\tens\id)
(\id\tens\id\tens\nabla)
(\check\sigma\tens\id\tens\id)(m\tens\nabla e) \cr
&=&  (\extd\tens\id\tens\id)
(\check\sigma\tens\id)(m\tens\nabla e)
-(\id\wedge\check\nabla\tens\id)
(\check\sigma\tens\id)(m\tens\nabla e) \cr
&&-\,(\id\wedge\check\sigma\tens\id)(\check\nabla\tens\id\tens\id)
(m\tens \nabla e)\cr
&&-\,(\id\wedge\check\sigma\tens\id)
(\check\sigma\tens\id\tens\id)(\id\tens\id\tens\nabla)(m\tens\nabla e)\ .
\end{eqnarray*}
Using property (b) of Definition~\ref{bhcsial} this becomes
\begin{eqnarray*}
\hat\nabla^{[1]}\hat\nabla(m\tens e)
&=&  (\extd\tens\id\tens\id)
(\check\sigma\tens\id)(m\tens\nabla e)
-(\id\wedge\check\nabla\tens\id)
(\check\sigma\tens\id)(m\tens\nabla e) \cr
&&\!\!\!\!-\,(\id\wedge\check\sigma\tens\id)(\check\nabla\tens\id\tens\id)
(m\tens \nabla e)\cr
&&\!\!\!\!-\,(\check\sigma\tens\id)(\id\tens\wedge\tens\id)
(\id\tens\id\tens\nabla)(m\tens\nabla e)\cr
&\!\!\!\!=&\!\!\!\!  (\extd\tens\id\tens\id)
(\check\sigma\tens\id)(m\tens\nabla e)
-(\id\wedge\check\nabla\tens\id)
(\check\sigma\tens\id)(m\tens\nabla e) \cr
&&\!\!\!\!-\,(\id\wedge\check\sigma\tens\id)(\check\nabla\tens\id\tens\id)
(m\tens \nabla e)
-(\check\sigma\tens\id)(m\tens (\id\wedge\nabla)\nabla e) \cr
&\!\!\!\!=&\!\!\!\!  (\extd\tens\id\tens\id)
(\check\sigma\tens\id)(m\tens\nabla e)
-(\id\wedge\check\nabla\tens\id)
(\check\sigma\tens\id)(m\tens\nabla e) \cr
&&\!\!\!\!-\,(\id\wedge\check\sigma\tens\id)(\check\nabla\tens\id\tens\id)
(m\tens \nabla e)
-(\check\sigma\tens\id) (\id\tens\extd \tens\id)(m\tens\nabla e)\ ,
\end{eqnarray*}
where we have used the flatness of $\nabla$ on $E$ in the last equality.
Now Lemma~\ref{lkvbihdz} completes the proof.\quad$\square$

\subsection{The bicategory of differentiable bimodules}\label{sec.bicat}
A possible way of understanding of  differentiable bimodules and 
induced functors between categories of connections is to construct  a 
suitable {\em bicategory}. Recall that a bicategory \cite{Ben:bic} 
consists of three layers of structures:  0-cells, 1-cells defined for 
any pair of 0-cells, and 2-cells defined for each pair of 1-cells. 
There are two types of composition: the horizontal composition of 
1-cells which is unital and associative up to isomorphisms and the 
vertical composition of 2-cells which is strictly associative and 
unital.   The following gathers all the data that constitute a 
bicategory relevant to differential bimodules.

\begin{defin}\label{def.bim}
The bicategory $\db$  of {\em differentiable bimodules} contains the 
following data:

(a) 0-cells are differential graded algebras $(\Omega^*A,\extd)$; we 
write $A$ for the zero-degree subalgebra of $\Omega^*A$.

(b) A 1-cell $\Omega^*A\to\Omega^*B$ is given by a differentiable 
bimodule $(M,\nabla_M,\sigma_M)$, i.e.\ $M$ is a $B$-$A$-bimodule, 
$\nabla_M : M\to \Omega^1B\otimes_B M$ is a left $B$-connection and 
$\sigma_M: M\otimes_A \Omega^1A\to \Omega^1B\otimes_B M$ is a 
generalised flip satisfying conditions of Definition~\ref{tenten}.

(c) A 2-cell $\xymatrix{\Omega^*A \ar[rr]^{(M,\nabla_M,\sigma_M)} 
&\ar@2[d]^ \phi & \Omega^*B \\
\Omega^*A \ar[rr]_{(N,\nabla_N,\sigma_N)} & & \Omega^*B}$ is given as 
a $B$-$A$ bimodule map $\phi: M\to N$ that commutes with covariant 
derivatives and generalised flip operators, i.e.\   such that 
$\nabla_N\circ \phi = (\id\otimes \phi)\circ\nabla_M$ and 
$\sigma_N\circ (\phi\otimes \id) = (\id\otimes \phi)\circ \sigma_M$.

\noindent The horizontal composition
$$
\xymatrix{\Omega^*A \ar[rr]^{(M,\nabla_M,\sigma_M)} && 
\Omega^*B\ar[rr]^{(N,\nabla_N,\sigma_N)} && \Omega^*C}
$$
is defined as a $C$-$A$ differentiable bimodule $(N\otimes_B M, 
\nabla_{N\otimes_B M},\sigma_{N\otimes_B M})$, where
$$
\nabla_{N\otimes_B M} = \nabla_N\otimes\id + (\sigma_N\otimes 
\id)\circ (\id\otimes\nabla_M), \quad \sigma_{N\otimes_BM} = 
(\sigma_N\otimes \id)\circ (\id\otimes\sigma_M).
$$
The vertical composition is the usual composition of mappings. The category of
1-cells $\Omega^*A\to\Omega^*B$ with morphisms provided by 2-cells is 
denoted by
$\db (\Omega^*A,\Omega^*B)$.
\end{defin}

It is left to the reader to check that the data collected in 
Definition~\ref{def.bim} indeed constitute a bicategory. Essentially 
this requires similar computations to those in the proof of 
Proposition~\ref{prop.functor}. The bicategory $\db$ contains all 
(left) connections in the following way.
\begin{lemma}\label{lemma.bim} View $k$ as a trivial differential 
graded algebra with the differential given by the zero map. Then
$$
\db (k,\Omega^*A) \equiv  {}_A\mathcal{E}.
$$
\end{lemma}
\noindent {\bf Proof:}\quad Since $\Omega ^1k =0$, every generalised 
flip $\sigma$ must be a zero map, thus an object in the category $\db 
(k,\Omega^*A)$ is a left $A$-module $M$ with a left
$A$-connection $\nabla_M: M\to \Omega^1A\otimes_A M$. As to the 
morphisms $\phi:M\to N$ in $\db (k,\Omega^*A)$, the commutativity 
with flips is trivially satisfied (as flips are zero maps), hence 
only the condition $\nabla_N\circ \phi = (\id\otimes 
\phi)\circ\nabla_M$ remains. This is equivalent to say that $\phi$ is 
a morphism in ${}_A\mathcal{E}$. \quad$\square$\medskip

In view of Lemma~\ref{lemma.bim}, the functor 
$(M,\nabla_M,\sigma_M)_*:  {}_A\mathcal{E}\to  {}_B\mathcal{E}$ 
constructed in Proposition~\ref{prop.functor} has a very simple and 
natural bicategorical explanation. Given 
a connection
$(E,\nabla_E) \in {}_A\mathcal{E}\equiv \db (k,\Omega^*A)$
and a differentiable bimodule
$(M,\nabla_M,\sigma_M)\in \db (\Omega^*A,\Omega^*B)$  one can 
construct a
differentiable bimodule in  $\db (k,\Omega^*B) \equiv 
{}_B\mathcal{E}$ as the horizontal composition of 1-cells
$$
\xymatrix{k \ar[rr]^{(E,\nabla_E)} && 
\Omega^*A\ar[rr]^{(M,\nabla_M,\sigma_M)} && \Omega^*B}.
$$
By the functoriality of the horizontal composition this results in a 
functor $ {}_A\mathcal{E}\to  {}_B\mathcal{E}$ described in 
Proposition~\ref{prop.functor}.

In a similar way one constructs a bicategory $\fdb$ of {\em flat 
differentiable bimodules} with the same 0-cells as in $\db$, the 
1-cells given as flat differentiable bimodules 
$(M,\nabla_M,\sigma^1_M,\sigma^2_M)$, where $\sigma^1_M$ and 
$\sigma^2_M$ are order one and two flip operators (cf.\ 
Definition~\ref{bhcsial}) and 2-cells $B$-$A$-bimodule maps commuting 
with $\nabla_M$, $\sigma^1_M$ and $\sigma^2_M$. The horizontal 
composition is given by
$$
\nabla_{N\otimes_B M} = \nabla_N\otimes\id + (\sigma^1_N\otimes 
\id)\circ (\id\otimes\nabla_M), \quad \sigma^i_{N\otimes_BM} = 
(\sigma^i_N\otimes \id)\circ (\id\otimes\sigma^i_M), \quad i=1,2,
$$
and the vertical composition is the usual composition of mappings. 
One easily shows that $ {}_A\mathcal{F}  \equiv \fdb (k,\Omega^*A)$ 
and then identifies the functor in Propostion~\ref{fourteen} as the 
horizontal composition of 1-cells in $\fdb$.
\subsection{The case of semi-free differential graded algebras.}
Recall that $\Omega^*A$ is said to be {\em semi-free} if and only if 
$\Omega^*A$ is isomorphic
to the tensor algebra of the $A$-bimodule $\Omega^1A$.
As observed in \cite{Roj:mat} there is a bijective correspondence 
between semi-free
differential graded algebras over $A$ and $A$-corings with a 
grouplike element (cf.\
\cite[29.8]{BrzWis:cor}). The constructions in Sections 2.1--2.3 have 
very natural
interpretation in terms of such corings and comodules. For more 
information on corings
and comodules we refer to \cite{BrzWis:cor}.

Starting with an $A$-coring $\cC$ and a grouplike element $g\in\cC$, we define
$\Omega^1 A = \ker\eC$, where $\eC:\cC\to A$ is the counit of $\cC$. 
The differential
is then defined by
     $ \extd(a) = ga - ag$, for all $a\in A$, and, for all $c^{1}\ut 
\cdots \ut c^{n}\in (\ker\eC)^{\otimes_A n}$,
     \begin{eqnarray*}
	\extd(c^{1}\ut \cdots \ut c^{n}) &=&
	g\ut c^{1}\ut \cdots \ut c^{n}+ (-1)^{n+1}c^{1}\ut \cdots \ut c^{n}
	\ut g \\
	&&+ \sum_{i=1}^{n}(-1)^{i}c^{1}\ut \cdots \ut
	c^{i-1}\ut \DC(c^{i})\ut c^{i+1}
	\ut \cdots\ut c^{n},
     \end{eqnarray*}
     where $\DC:\cC\to\cC\ot_A\cC$ is the coproduct in $\cC$. The {\em 
density condition} for $\Omega^1A$, i.e.\ the requirement that any 
one-form is a linear
     combination of $a\extd a'$, is equivalent to the requirement that 
the map $A\ot A\to \cC$, $a\ot a'\mapsto aga'$ be surjective (note 
the similarity with the definition of
     a {\em space cover}  in \cite{KonRos:smo}).

Let $E$ be a left $A$-module.    As explained in
     \cite[29.11]{BrzWis:cor}, connections $\nabla: E\to 
\Omega^1A\ot_A E = \ker\eC\ot_A E$  are in bijective correspondence 
with left $A$-module sections
     of $\eC\ot \id:\cC\ot_A E\to E$, i.e.\ left $A$-linear maps 
$\varrho^E: E\to \cC\ot_A E$
     such that $(\eC\ot\id)\circ \varrho^E = \id$. Furthermore, flat 
connections are
     in bijective correspondence with left $\cC$-coactions in $E$. 
This correspondence, explicitly given by
     $$
     \varrho^E(e) = g\ot e - \nabla(e), \qquad \forall e\in E,
     $$
     establishes an isomorphism of the categories of flat connections on $A$ and
     left $\cC$-comodules.

Let $\cC$ be an $A$-coring with a grouplike element $g_\cC$ and $\cD$ be
a $B$-coring with a grouplike element $g_\cD$. Recall that a morphism 
of corings
consists of an algebra map $\theta_0:A\to B$ and an $A$-bimodule map 
$\theta_1: \cC\to \cD$ that respects the coproducts and counits (cf.\ 
\cite[24.1]{BrzWis:cor}
for more details). Any morphism of corings $(\theta_0,\theta_1)$ such 
that $\theta_1(g_\cC) = g_\cD$ is a differentiable algebra map. 
Incidentally, such a morphism of 
corings is termed a {\em morphism of space covers} in 
\cite{KonRos:smo}. Let $\nabla: E\to \Omega^1A\ot_A E$ be a 
connection, and
$\varrho^E: E\to \cC\ot_A E$ be the corresponding section  of 
$\eC\ot\id$. Then the
section $\varrho^{B\ot_AE}: E\to \cD\ot_B E$ of $\eD\ot\id$ 
corresponding to the induced connection in $B\ot_A E$ comes out as,
$$
\varrho^{B\ot_AE}(b\ot e) =  b\theta_1(e\sco{-1})\ot e\sco 0,
$$
where $\varrho^E(e) =  e\sco{-1}\ot e\sco 0$ (summation implicitly 
understood). In view of the isomorphism
${}_A\mathcal{F} \cong {}^\cC\M$, the corresponding functor between 
the categories of flat connections  described in Theorem~\ref{hjhjhj} 
can be identified with the induction functor between categories of 
left comodules (cf.\ \cite[24.6]{BrzWis:cor}).

For differential graded algebras corresponding to an $A$-coring $\cC$ 
with a grouplike element $g_\cC$ and a $B$-coring $\cD$  with a 
grouplike element $g_\cD$, differentiable bimodules $(M,{\nabla}, 
{\sigma})$ are in bijective
correspondence with pairs $(M,\Phi)$, where $M$ is a $(B,A)$-bimodule 
and $\Phi: M\ot_A\cC\to \cD\ot_B M$ is a $(B,A)$-bimodule map 
rendering commutative the following diagram
\begin{equation}\label{diag.eps}
\xymatrix{ M\otimes_A\cC\ar[rr]^\Phi \ar[dr]_{\id\ot\eC} && \cD\ot_B 
M\ar[dl]^{\eD\ot \id} \\
& M&}
\end{equation}
Furthermore, differentiable flat bimodules $(M,{\nabla}, {\sigma})$ 
are in bijective
correspondence with pairs $(M,\Phi)$ such that in addition to 
(\ref{diag.eps}) also
the following diagram
\begin{equation}\label{diag.del}
\xymatrix{ M\otimes_A\cC\ar[rr]^\Phi \ar[d]_{\id\ot\DC} && \cD\ot_B 
M\ar[d]^{\DD\ot \id} \\
M\ot_A\cC\ot_A\cC\ar[dr]_{\Phi\ot\id}&& \cD\ot_B\cD\ot_B M \\
&\cD\ot_B M\ot_A\cC\ar[ru]_{\cD\ot\Phi}&}
\end{equation}
is commutative. The correspondence is given by $\sigma = 
\Phi\mid_{M\ot_A\ker\eC}$ and
$$
\Phi (m\ot c) = g_\cD\ot m\eC(c) - \nabla(m)\eC(c) +\sigma(m\ot 
(c-g_\cC\eC(c))),
$$
for all  $m\in M$ and $c\in \cC$. An interesting point to note here 
is that the map $\Phi$ is well defined, i.e.\ factors through the 
coequaliser defining $M\ot_A\cC$, thanks to the last condition in 
Definition~\ref{tenten} (the compatibility between connection and 
$\sigma$).

A pair $(M,\Phi)$ satisfying conditions (\ref{diag.eps}), 
(\ref{diag.del}) constitutes a 1-cell in the {\em left bicategory of 
corings} $\LEM$  defined in \cite{BrzElK:bic} as the
  bicategory of comonads in the bicategory $\mathsf{Bim}$ of rings and 
bimodules following general procedure in \cite{Str:for}, 
\cite{LacStr:for}.
  In view of the discussion of Section~\ref{sec.bicat} and the present 
section,  $\LEM$ can be understood as a full sub-bicategory of $\fdb$.

\section{The long exact sequence}
Consider a short exact sequence $0\to E\to F\to G\to 0$ in
${}_A\mathcal{F}$,
and suppose that the modules
$\Omega^* A$ are flat
(i.e.\ tensoring with them preserves exactness).
We assume these conditions for the remainder of the section.
  From this we form the following diagram,
where the rows are exact, and the columns form cochain complexes
(i.e.\ the vertical maps compose to give zero).

\begin{picture}(270,140)(-45,-72)

\put(0,47){$0$}
\put(8,50){\vector(1,0){38}}
\put(48,47){$E$}
\put(60,50){\vector(1,0){70}}
\put(132,47){$F$}
\put(144,50){\vector(1,0){70}}
\put(216,47){$G$}
\put(228,50){\vector(1,0){38}}
\put(270,47){$0$}
\put(88,53){$\phi$}
\put(174,53){$\psi$}

\put(52,42){\vector(0,-1){25}}
\put(135,42){\vector(0,-1){25}}
\put(219,42){\vector(0,-1){25}}
\put(57,28){$\nabla$}
\put(140,28){$\nabla$}
\put(224,28){$\nabla$}

\put(0,7){$0$}
\put(8,10){\vector(1,0){13}}
\put(25,7){$\Omega^1 A\tens_A E$}
\put(77,10){\vector(1,0){27}}
\put(108,7){$\Omega^1 A\tens_A F$}
\put(158,10){\vector(1,0){29}}
\put(192,7){$\Omega^1 A\tens_A G$}
\put(242,10){\vector(1,0){23}}
\put(270,7){$0$}
\put(77,13){$\id\tens\phi$}
\put(159,13){$\id\tens\psi$}

\put(52,2){\vector(0,-1){25}}
\put(135,2){\vector(0,-1){25}}
\put(219,2){\vector(0,-1){25}}
\put(57,-13){$\nabla^{[1]}$}
\put(140,-13){$\nabla^{[1]}$}
\put(224,-13){$\nabla^{[1]}$}

\put(0,-33){$0$}
\put(8,-30){\vector(1,0){13}}
\put(25,-33){$\Omega^2 A\tens_A E$}
\put(77,-30){\vector(1,0){27}}
\put(108,-33){$\Omega^2 A\tens_A F$}
\put(158,-30){\vector(1,0){29}}
\put(192,-33){$\Omega^2 A\tens_A G$}
\put(242,-30){\vector(1,0){23}}
\put(270,-33){$0$}
\put(77,-27){$\id\tens\phi$}
\put(159,-27){$\id\tens\psi$}

\put(52,-38){\vector(0,-1){25}}
\put(135,-38){\vector(0,-1){25}}
\put(219,-38){\vector(0,-1){25}}
\put(57,-53){$\nabla^{[2]}$}
\put(140,-53){$\nabla^{[2]}$}
\put(224,-53){$\nabla^{[2]}$}

\end{picture}

What follows is standard homological algebra, but not all readers may 
be familiar with it. Note that for (e.g.) $\psi:F\to G$
  we write  $\psi^{-1}(g)$ for $g\in G$
to mean a choice of $f\in F$ for which $\psi(f)=g$. It will turn out that
the maps eventually defined by using this
such potentially multivalued maps will turn out to be unique,
and we have no wish to introduce the complication of topologised
cochain complexes, and so have no need to worry about the continuity
of the resulting operations. It is merely notation used to try to clarify the
definitions and proofs. Again take $\Gamma E=\{e\in E:\nabla e=0\}$.

\begin{propos}
The sequence $0\to \Gamma E\to \Gamma F\to\Gamma G$ is exact.
\end{propos}
\noindent {\bf Proof:}\quad
It is immediate that
$\phi:\Gamma E\to \Gamma F$ is one-to-one, and that the composition
$\Gamma E\to \Gamma F\to \Gamma G$ is zero. To show that
$\Gamma E\to \Gamma F\to \Gamma G$ is exact, take $f\in \Gamma F$
with $\psi (f)=0$. As $E\to F\to G$ is exact, there is an $e\in E$ with
$\phi (e)=f$. By following the top left
commutative square in the diagram and using the fact that
$\id\tens\phi:\Omega^1 A\tens_A E\to \Omega^1 A\tens_A F$
is one-to-one we see that $\nabla e=0$. \quad $\square$

\begin{propos}
The (multivalued) map $(\id\tens\phi)^{-1}\nabla\psi^{-1}:
\Gamma G\to \Omega^1 A\tens_A E$ quotients to a well defined
connecting map $\Gamma G\to H^1(A;E)$.
\end{propos}
\noindent {\bf Proof:}\quad Begin with $g\in\Gamma G$, and take an $f\in
F$ with $\psi(f)=g$. By using the top right commutative square in the
diagram, $\nabla f\in\ker(\id\tens\psi:\Omega^1 A\tens_A F\to
\Omega^1 A\tens_A G)$.
Then by the exactness of the rows, there is an $x\in \Omega^1 A\tens_A E$
with $(\id\tens \phi)(x)=\nabla f$. By exactness of the second row,
to show that $e'\in\ker\nabla^{[1]}$ we only have to show that
$\nabla^{[1]}(\id\tens\phi)(x)=0$, i.e.\ that $\nabla^{[1]}\nabla f=0$,
which is true. Then $[x]\in H^1(A;E)$, but now we ask if it is unique.

Suppose that we have $f'\in F$  with $\psi(f')=g$, and
$x'\in \Omega^1 A\tens_A E$
with $(\id\tens \phi)(x')=\nabla f'$. Then $f'-f=\phi(e)$ for some $e\in E$,
and $(\id\tens \phi)(x'-x)=\nabla 
(f'-f)=\nabla\phi(e)=(\id\tens\phi)(\nabla e)$.
As $\id\tens\phi$ is one-to-one we deduce that $x'-x=\nabla e$. \quad $\square$

\begin{remark}
As this is not a text on homological algebra, we will now merely quote the
result of continuing with the methods outlined: Given the conditions 
at the beginning of this section, there is a long exact sequence
\begin{eqnarray*}
H^0(A,E)\tto H^0(A,F)\tto\!\!\!\!\!\!\!\!&& H^0(A,G)\tto 
H^1(A,E)\\&&\tto H^1(A,F)
\tto H^1(A,G)
\tto
H^2(A,E)\tto  \dots
\end{eqnarray*}

\end{remark}

\section{Noncommutative fibre bundles}
We consider a possible meaning for a
differentiable algebra map $\iota:B\to X$ to be a `fibration'
with `base algebra' $B$ and `total algebra' $X$.
From here we will 
require that the differential calculi satisfy the density condition.

\begin{defin} \label{hhhddd}
Define the cochain complexes
\begin{eqnarray*}
\Xi_m^0 X\ =\ \iota_*\Omega^m B . X\ ,\quad
\Xi_m^n X\ =\ \frac
{\iota_*\Omega^m B \wedge \Omega^n X}
{\iota_*\Omega^{m+1} B \wedge \Omega^{n-1} X} \quad ( n>0)\ ,
\end{eqnarray*}
with differential $\extd:\Xi_m^n X\to \Xi_m^{n+1} X$ defined by
$\extd[\omega]_m=[\extd \omega]_m$, where $\omega\in
\iota_*\Omega^m B \wedge \Omega^n X$ and $[~]_m$ is the
corresponding quotient map.

The maps $\Theta_m:\Omega^m B\tens_B \Xi_0^n X \to \Xi_m^{n} X$
defined by $\Theta_m(\omega\tens[\xi]_0)=[\iota_*\omega\wedge\xi]_m$
are cochain maps if $\Omega^m B\tens_B \Xi_0^* X$
is given the differential $(-1)^m\id\tens\extd$.
\end{defin}

\begin{remark}
  To see that the differential in Definition~\ref{hhhddd}
  is well defined, note that for all $m,n\ge 0$,
  $\extd$ maps $\iota_*\Omega^m B \wedge \Omega^n X$ into
  $\iota_*\Omega^m B \wedge \Omega^{n+1} X$. This is because
  $\extd \Omega^m B\subset \Omega^{m+1} B\subset
  \Omega^m B\wedge\Omega^1 B$ (note the use of the density condition here).

There is a left $B$-module structure for $\Xi_m^n X$ given by
$b.\xi=\iota(b)\xi$.
  As $\extd (\iota(b).\theta)=\iota_*(\extd b)\wedge \theta+
\iota(b).\extd\theta$, we see that $\extd:\Xi_m^n X\to \Xi_m^{n+1} X$
is a left $B$-bimodule map, so the cohomology $H^n(\Xi_m^* X)$
inherits a left $B$-bimodule structure.
\end{remark}

In this degree of generality, this construction might be merely curious,
but consider an example:

\begin{example}
Let $X=B\tens F$ where $F$ is an algebra with differential structure,
  and give $X$ the tensor product differential structure. By definition
  $\iota(b)=b\tens 1$ and
  \begin{eqnarray*}
\Omega^n X \,=\,(\Omega^0 B\tens\Omega^n F) \oplus\dots\oplus
(\Omega^n B\tens\Omega^0 F)\ ,
\end{eqnarray*}
so there is an isomorphism of cochain complexes
$B\tens\Omega^n F\to \Xi^n_0 X$ given by $b\tens\xi\mapsto \iota(b)\xi$. It
follows that $H^n(\Xi^*_0 X)$ is just $B\tens H_{dR}^n(F)$,
the fibre cohomology module. Also
this module has a flat left $B$-connection
$\nabla: B\tens H^*_{dR}(F)\to \Omega^1 B \tens_B B\tens H^*_{dR}(F)$
given by
$\nabla(b\tens x)=\extd b\tens 1\tens x$. The de Rham cohomology of
$B$ with coefficients in this module with flat connection
is $H^*_{dR}(B)\tens H^*_{dR}(F)$, which by the K\"unneth theorem is
just the cohomology of $X=B\tens F$.
\end{example}

In topology fibrations can be built from open covers of the base space,
and a trivial fibration over each open set. Our example has just dealt with
what would be a noncommutative trivial fibration, so we might ask what
a more general noncommutative fibration would look like.
By analogy we might consider $\Xi^n_0 X$ to be the `vertical' or `fibre'
forms, and its cohomology to be the cohomology of the `fibre'
of the map. In the topological case, this cohomology can form a
non-trivial bundle
over the base space. We have seen that for noncommutative de Rham
cohomology it is reasonable to have coefficient bundles with flat connection,
and this is the route that we will take for our version of a fibration.

\begin{propos}\label{uioy}
  Suppose that $\Theta_1:\Omega^1 B\tens_B \Xi_0^* X \to
\Xi_1^{*} X$ (as defined in Definition~\ref{hhhddd}) is invertible. Then there 
is a left $B$-covariant derivative
$\nabla: H^n(\Xi_0^* X)\to \Omega^1 B\tens_B H^n(\Xi_0^* X)$
defined by $[\omega] \mapsto (\id\tens [~])\Theta_1^{-1}
[\extd\omega]_1$.
\end{propos}
\noindent {\bf Proof:}\quad
If $[\omega]_0\in Z^n=\ker\extd:\Xi_0^n X \to \Xi_0^{n+1} X$, then
$\extd\omega\in \iota_*\Omega^1 B \wedge \Omega^n X$.
  Then $[\extd\omega]_1\in\Xi_1^{n} X$ is a cocycle,
so $(\id\tens \extd)\Theta_1^{-1}[\extd \omega]_1=0\in
\Omega^1 B\tens_B \Xi_0^{n+1} X$,
i.e.\ $\Theta_1^{-1}[\extd \omega]_1\in\Omega^1 B\tens_B Z^n$.

Now suppose that $[\omega']=[\omega]\in Z^n$. Then
$\omega'-\omega\in \iota_*\Omega^1 B \wedge \Omega^{n-1} X$,
so we get $\Theta_1^{-1}[\omega'-\omega]_1\in
\Omega^1 B\tens_B \Xi_0^{n-1} X$. As $\Theta^{-1}$ is a cochain map,
$-(\id\tens\extd)\Theta_1^{-1}[\omega'-\omega]_1=
\Theta_1^{-1}[\extd\omega'-\extd\omega]_1=
\Theta_1^{-1}[\extd\omega']_1-\Theta_1^{-1}[\extd\omega]_1$.
Thus $\Theta_1^{-1}[\extd\omega']_1-\Theta_1^{-1}[\extd\omega]_1
\in \Omega^1 B\tens_B \extd\Xi_0^{n-1} X$, so we get a well
defined map $Z^n\to \Omega^1 B\tens_B H^n(\Xi_0^* X)$.

To complete showing that
$\nabla$ is well defined, we show that $\extd\Xi_0^{n} X$ maps to zero, which
we see as $\nabla[\extd\xi]=
(\id\tens [~])\Theta_1^{-1}[\extd^2\xi]_1=0$.

Finally we need to show the left connection condition:
\begin{eqnarray*}
\nabla[\iota(b).\omega] &=& (\id\tens [~])\Theta_1^{-1}[\iota(\extd 
b)\wedge \omega+
\iota(b).\extd \omega]\,=\,\iota(b)\wedge[\omega]+b.\nabla[\omega]\ 
.\quad\square
\end{eqnarray*}

\begin{propos}
Suppose that $\Theta_m:\Omega^m B\tens_B \Xi_0^* X \to
\Xi_m^{*} X$ (as defined in Definition~\ref{hhhddd}) is invertible for $m=1,2$. Then
the curvature of the connection on $H^n(\Xi_0^* X)$ described in
Proposition~\ref{uioy} is zero.
\end{propos}
\noindent {\bf Proof:}\quad Take
$[\omega]_0\in Z^n=\ker\extd:\Xi_0^n X \to \Xi_0^{n+1} X$, and write
$\Theta_1^{-1}[\extd\omega]_1=\sum_i \xi_i\tens[\eta_i]_0\in
\Omega^1 B\tens_B Z^n$. Likewise write
$\Theta_1^{-1}[\extd\eta_i]_1=\sum_j \chi_{ij}\tens[\mu_{ij}]_0\in
\Omega^1 B\tens_B Z^n$. Now write the composition $\nabla^{[1]}
\nabla$ as
\begin{eqnarray*}
[\omega]\mapsto \sum_i \xi_i\tens[\eta_i] \mapsto
\sum_i\Big(\extd\xi_i\tens[\eta_i] -\sum_j
\xi_i\wedge \chi_{ij} \tens [\mu_{ij}]\Big)\ .
\end{eqnarray*}
Inserting the definition of $\Theta^{-1}$, we get
$[\extd\omega]_1=\sum_i[\iota_*\xi_i\wedge \eta_i]_1$
and
$[\extd\eta_i]_1=\sum_j[\iota_*\chi_{ij}\wedge \mu_{ij}]_1$.
This means that $\extd\omega-\sum_i\iota_*\xi_i\wedge \eta_i
\in \iota_*\Omega^2 B\wedge \Omega^{n-1}X$, so we write
$\extd\omega-\sum_i\iota_*\xi_i\wedge \eta_i=
\sum_k \iota_*\tau_k\wedge\lambda_k$ where $\tau_k\in\Omega^2 B$
and $\lambda_k\in \Omega^{n-1}X$. Applying $\extd$ to this, we get
\begin{eqnarray*}
\sum_k \Big(\iota_*\extd\tau_k\wedge\lambda_k+
  \iota_*\tau_k\wedge\extd\lambda_k\Big) &=&
  \sum_i\Big(\iota_*\xi_i\wedge \extd\eta_i-\iota_*\extd\xi_i\wedge 
\eta_i\Big)\ .
\end{eqnarray*}
Then we obtain $\sum_k[\iota_*\tau_k\wedge\extd\lambda_k]_2=
\sum_{ij}[\iota_*(\xi_i\wedge\chi_{ij})\wedge\mu_{ij}]_2
-\sum_i[\iota_*(\extd\xi_i)\wedge\eta_i]_2$. Then the two elements
$ \sum_k\tau_k\tens[\extd\lambda_k]_0 $ and
$ \sum_{ij}\xi_i\wedge\chi_{ij}\tens[\mu_{ij}]_0
-\sum_i\extd\xi_i\tens[\eta_i]_0 $
of $\Omega^2 B\tens_B \Xi_0^{n-1}$ map to the same thing under $\Theta_2$,
so by our assumption they must be equal. Now as $[\extd\lambda_k]_0$ is
a coboundary, the curvature must vanish. \quad$\square$

\section{Spectral sequences}
The reader should refer to \cite{spseq}  for the details of the homological
algebra used to construct the spectral sequence. We will
merely quote the results.

\begin{remark} \label{sprem}
Start with a differential graded module $C^n$ ( for $n\ge 0$) and 
$\extd :C^n\to
C^{n+1}$ with $\extd^2=0$. Suppose that $C$ has a filtration
$F^m C\subset C=\oplus_{n\ge 0}C^n$ for
$m\ge 0$ so that:

(1)\quad $\extd F^m C \subset F^m C$ for all $m\ge 0$
(i.e.\ the filtration is preserved by $\extd$);

(2)\quad $F^{m+1} C\subset F^m C$ for all $m\ge 0$
(i.e.\ the filtration is decreasing);

(3)\quad $F^0 C=C$ and $F^m C^n=F^m C\cap C^n=\{0\}$
for all $m>n$ (a boundedness condition).

\noindent
Then there is a spectral sequence $(E_r^{*,*}, \extd_r)$
for $r\ge 1$ with $\extd_r$ of bidegree $(r,1-r)$ and
$$
E_1^{p,q}\,=\,H^{p+q}(F^pC/F^{p+1}C)\,=\,
\frac{{\rm ker}\,\extd:F^pC^{p+q}/F^{p+1}C^{p+q}\to 
F^pC^{p+q+1}/F^{p+1}C^{p+q+1}}
{{\rm im}\,\extd:F^pC^{p+q-1}/F^{p+1}C^{p+q-1}\to F^pC^{p+q}/F^{p+1}C^{p+q}}\ .
$$
In more detail, we define
\begin{eqnarray*}
Z_r^{p,q} &=& F^{p} C^{p+q} \cap \extd^{-1}(F^{p+r} C^{p+q+1})\ ,\cr
B_r^{p,q} &=& F^{p} C^{p+q} \cap \extd(F^{p-r} C^{p+q-1})\ ,\cr
E_{r}^{p,q} &=& Z_r^{p,q}/(Z_{r-1}^{p+1,q-1}+B_{r-1}^{p,q})\ .
\end{eqnarray*}
The differential $\extd_r:E_{r}^{p,q}\to E_{r}^{p+r,q-r+1}$ is the
map induced on quotienting $\extd:Z_{r}^{p,q}\to Z_{r}^{p+r,q-r+1}$.

The spectral sequence converges to $H^*(C,\extd)$ in the sense that
\begin{eqnarray*}
E_\infty^{p,q} \cong \frac
{F^pH^{p+q}(C,\extd)}
{F^{p+1}H^{p+q}(C,\extd)}\ ,
\end{eqnarray*}
where $F^pH^*(C,\extd)$ is the image of the
map $H^*(F^p C,\extd)\to H^*(C,\extd)$
induced by inclusion $F^p C\to C$.
\end{remark}

Now take the case of a differentiable algebra map
$\iota:B\to X$. We can give the following example of a spectral sequence.

\begin{remark} \label{deffil}
Define the filtration $F^m \Omega^{n+m} X
=\iota_*\Omega^m B \wedge \Omega^n X$
of $\Omega^* X$. This obeys conditions (1) and (2) of Remark~\ref{sprem}
as\begin{eqnarray*}
\iota_*\Omega^{m+1} B \wedge \Omega^n X \subset
\iota_*\Omega^{m} B \wedge \iota_*\Omega^{1} B \wedge \Omega^n X \subset
\iota_*\Omega^{m} B \wedge \Omega^{n+1} X \ .
\end{eqnarray*}
We have boundedness as $\iota_*\Omega^{0} B \wedge \Omega^n X=
\Omega^n X$, and by convention $\Omega^{n}X=0$ for $n<0$.
Note that
\begin{eqnarray*}
\frac{F^p\Omega^{p+q}X}{F^{p+1}\Omega^{p+q}X} \,=\,\Xi_p^q X\ ,
\end{eqnarray*}
and we obtain a spectral sequence with $E_1^{p,q}\cong H^q(\Xi_p^* X)$
which converges to $H^*_{dR}(X)$ in the sense described in Remark~\ref{sprem}.
The differential $\extd_1:H^q(\Xi_p^* X)\to H^q(\Xi_{p+1}^* X)$ is the map
given by applying $\extd$ to cocycles in $\Xi_p^* X$,
taking care what space you end up in!
\end{remark}

\begin{defin} \label{jasdfvajkhv}
The differentiable algebra map $\iota:B\to X$ is called a {\em 
differential fibration}
if $\Theta_m:\Omega^m B\tens_B \Xi_0^* X \to
\Xi_m^{*} X$ (as given in Definition~\ref{hhhddd})  is invertible for all $m\ge 0$.
\end{defin}

\begin{theorem}
Suppose that $\iota:B\to X$ is a differential fibration. Then there 
is a spectral sequence converging to $H^*_{dR}(X)$ with
\begin{eqnarray*}
E_2^{p,q} \cong H^p(B;H^q(\Xi_0^* X),\nabla)
\end{eqnarray*}
\end{theorem}
\noindent {\bf Proof:}\quad We note that
$\Theta_{m*}:\Omega^p B\tens_B
H^q(\Xi_0^* X) \to H^q(\Xi_p^* X)$ is an isomorphism, and that it commutes
with the differential in the spectral sequence if we use the flat
connection cochain complex on $\Omega^p B\tens_B H^q(\Xi_0^* X)$.
\quad$\square$

\section{The multiplicative structure}
Even if one is not {\it a priori}  interested in a multiplicative
structure on the cohomology theories, in algebraic topology
a knowledge of the multiplicative structure can help find the
differentials in the spectral sequence.
In this section we suppose that the differentiable algebra
  map $\iota:B\to X$ is a differential fibration,
and that the following condition holds:

\begin{defin} \label{dbcccc}
The map $\iota:B\to X$ will be said to satisfy the {\em differential 
braiding condition} if
$ \Omega^n X\wedge \iota_*\Omega^m B\subset \iota_*\Omega^m B \wedge 
\Omega^n X$
for all $n,m\ge 0$.
\end{defin}

\begin{remark}
Note that the condition in Definition~\ref{dbcccc} means that the wedge 
multiplication
preserves the filtration in the construction of the spectral sequence, as
\begin{eqnarray*}
( \iota_*\Omega^i B \wedge \Omega^j X)\wedge
( \iota_*\Omega^k B \wedge \Omega^l X)  &\subset&
\iota_* \Omega^i B \wedge\iota_*\Omega^k B \wedge
\Omega^j X\wedge \Omega^l X\cr
&\subset & \iota_*\Omega^{i+k} B \wedge \Omega^{j+l} X\ ,
\end{eqnarray*}
so there is a multiplicative structure on the spectral sequence.
However we have gone to considerable trouble to show that the $E_2$ page of the
spectral sequence can be expressed in terms of a cohomology bundle
with connection, so we shall look at what this multiplicative structure means
in these terms.
\end{remark}

\begin{propos}
Define a map $\hat\sigma:\Xi^n_0 X\tens_B \Omega^m B\to \Omega^m B\tens_B
\Xi^n_0 X$ by $\hat\sigma([\xi]_0\tens \omega)=\omega'_i\tens[\xi'_i]_0$
(summation implicit), where $[\iota_*\omega'_i\wedge \xi'_i]_m=(-1)^{nm}\,[\xi
\wedge\iota_*\omega]_m$. For the cochain structure on $\Xi^*_0$, 
$$
\hat\sigma((\ker \extd)\tens_B \Omega^m B)\subset \Omega^m 
B\tens_B(\ker \extd) \quad
\textrm{and} \quad
\hat\sigma(({\rm im}\, \extd)\tens_B \Omega^m B)
\subset \Omega^m B\tens_B({\rm im}\, \extd), 
$$
so there is a well defined map
$\sigma:H^n(\Xi^*_0 X)\tens_B \Omega^m B\to \Omega^m B\tens_B
H^n(\Xi^*_0)$.
\end{propos}
\noindent {\bf Proof:}\quad First suppose that $[\xi]_0\in \ker 
\extd\subset \Xi^n_0 X$.
We write $\hat\sigma([\xi]_0\tens \omega)=\omega'_i\tens[\xi'_i]_0$, where
\begin{eqnarray*}
\iota_*\omega'_i\wedge \xi'_i\cong(-1)^{nm}\,\xi\wedge\iota_*\omega
\ \; {\rm modulo}\ \;  \iota_*\Omega^{m+1} B\tens_B \Omega^{n-1} X\ .
\end{eqnarray*}
On applying $\extd$,
\begin{eqnarray*}
\iota_*\extd\omega'_i\wedge \xi'_i+(-1)^m\,
\iota_*\omega'_i\wedge \extd\xi'_i&\cong& 
(-1)^{nm}\,\extd\xi\wedge\iota_*\omega
+(-1)^{nm+n}\,\xi\wedge\iota_*\extd\omega\\
&&~~~~~~~~~~~~~~~~~~~~\
  {\rm modulo}\ \;\iota_*\Omega^{m+1} B\tens_B \Omega^{n} X\ .
\end{eqnarray*}
As $\extd\xi\in\iota_*\Omega^1 B\tens_B \Omega^n X$, this shows that
$[\iota_*\omega'_i\wedge \extd\xi'_i]_m=0$, therefore, the fibration condition gives
$\omega'_i\tens_B [\extd\xi'_i]_0=0$.

Now take $[\eta]_0\in \Xi_0^{n-1}X$, and then  find
\begin{eqnarray*}
\iota_*\omega'_i\wedge \eta'_i\cong \eta\wedge\iota_*\omega
\ \;{\rm modulo}\ \;\iota_*\Omega^{m+1} B\tens_B \Omega^{n-2} X\ .
\end{eqnarray*}
Applying $\extd$ gives
\begin{eqnarray*}
\iota_*\extd\omega'_i\wedge \eta'_i+(-1)^m\,\iota_*\omega'_i\wedge \extd\eta'_i
&\cong& \extd\eta\wedge\iota_*\omega -  (-1)^n\,\eta\wedge\iota_*\extd\omega
\\
&& \qquad {\rm modulo}\;\; \iota_*\Omega^{m+1} B\tens_B \Omega^{n-1} X\ ,
\end{eqnarray*}
which reduces to
\begin{eqnarray*}
(-1)^m\,\iota_*\omega'_i\wedge \extd\eta'_i
\cong \extd\eta\wedge\iota_*\omega
\ {\rm modulo}\ \iota_*\Omega^{m+1} B\tens_B \Omega^{n-1} X\ .\quad\square
\end{eqnarray*}

\begin{propos}
If the differential braiding condition holds, then there exists a well defined
map $\wedge:\Xi^r_0 X\tens \Xi^s_0 X\to \Xi^{r+s}_0 X$ defined by
$[\xi]_0\wedge[\eta]_0=[\xi\wedge\eta]_0$, and this gives a well
defined map $\wedge:H^r(\Xi^*_0 X)\tens H^s(\Xi^*_0 X)\to H^{r+s}(\Xi^*_0 X)$.
\end{propos}
\noindent {\bf Proof:}\quad To show that the map $[\xi]_0\tens [\eta]_0\mapsto
[\xi\wedge\eta]_0$ is well defined we need to show that both
$\iota_* \Omega^1 B\wedge\Omega^{r-1} X\wedge\Omega^s X$ and
$\Omega^{r} X\wedge\iota_* \Omega^1 B\wedge\Omega^{s-1} X$
are contained in $\iota_* \Omega^1 B\wedge\Omega^{r+s-1} X$. The 
first inclusion
is automatic, and the second follows from the differential braiding condition.
The rest is left to the reader. \quad$\square$

\begin{propos}
For all $x\in H^n(\Xi^*_0 X)$ and $\omega\in \Omega^m B$,
\begin{eqnarray*}
(\wedge\tens\id)(\id\tens\sigma)(\nabla x\tens\omega) 
+\sigma(x\tens\extd\omega)
\,=\,[\extd\tens\id+(-1)^m(\id\wedge\nabla)]\sigma(x\tens\omega)\ .
\end{eqnarray*}
\end{propos}
\noindent {\bf Proof:}\quad
For all $\omega\in\Omega^m B$ and $\xi\in\Omega^n X$, we have defined
$\hat\sigma(\xi\tens\omega)=\omega_i\tens\xi_i$, where
\begin{eqnarray*}
(-1)^{nm}\, \xi\wedge\iota_*\omega\,=\,\iota_*\omega_i\wedge 
\xi_i+\iota_*\phi_i\wedge
\eta_i\ ,
\end{eqnarray*}
for some $\phi_i\in\Omega^{m+1}B$ and $\eta_i\in\Omega^{n-1}X$.
Taking $\extd$ of this gives
\begin{eqnarray}\label{hyre}
(-1)^{nm}\, \extd\xi\wedge\iota_*\omega + (-1)^{nm+n}\, 
\xi\wedge\iota_*\extd\omega
  &=& \iota_*\extd\omega_i\wedge \xi_i+(-1)^m\,\iota_*\omega_i\wedge 
\extd\xi_i \cr
&& 
+\,\iota_*\extd\phi_i\wedge\eta_i+(-1)^{m+1}\,\iota_*\phi_i\wedge\extd\eta_i\ 
.
\end{eqnarray}
Now we suppose that $[\xi]_0\in\ker\extd:\Xi^n_0 X\to \Xi^{n+1}_0 X$, and then
we also have $[\extd\xi_i]_m=0$. This means that all the terms of (\ref{hyre})
are in $\iota_*\Omega^{m+1}B\wedge \Omega^{n}X$,
and using the quotient map $[.]_{m+1}$ we obtain
\begin{eqnarray}\label{hyre2}
(-1)^{nm}\, [\extd\xi\wedge\iota_*\omega]_{m+1} \!\!\!\!\!\!\!\!&&+\ 
(-1)^{nm+n}\, [\xi\wedge\iota_*\extd\omega]_{m+1}
  = [\iota_*\extd\omega_i\wedge \xi_i]_{m+1}\cr
  &&+\
  (-1)^m\,[\iota_*\omega_i\wedge \extd\xi_i ]_{m+1}
  +\,(-1)^{m+1}\,[\iota_*\phi_i\wedge\extd\eta_i]_{m+1}\ .
\end{eqnarray}
Now write $\nabla \xi=\psi_i\tens[\zeta_i]_0\in\Omega^1B\tens\Xi^n_0$ and
$\nabla \xi_i=\psi_{ik}\tens[\zeta_{ik}]_0$.  Substituting this in 
(\ref{hyre2}) gives
\begin{eqnarray}\label{hyre3}
(-1)^{nm}\, [\iota_*\psi_i\wedge\zeta_i\wedge\iota_*\omega]_{m+1}
  &=& [\iota_*\extd\omega_i\wedge \xi_i]_{m+1}+
  (-1)^m\,[\iota_*\omega_i\wedge \iota_*\psi_{ik}\wedge\zeta_{ik} ]_{m+1}\cr
   &&\!\!\!\!\!\!\!\!\!\!\!\!\!\!+ (-1)^{nm+n}\, 
[\xi\wedge\iota_*\extd\omega]_{m+1}
  +\,(-1)^{m+1}\,[\iota_*\phi_i\wedge\extd\eta_i]_{m+1}.
\end{eqnarray}
On passing to the cohomology the last term in (\ref{hyre3}) vanishes,
giving the result.\quad$\square$

\begin{propos}
For $x,y\in H^*(\Xi^*_0 X)$, $\nabla(x\wedge y)=\nabla x\wedge y
+(\sigma\wedge\id)(x\tens\nabla y).$
\end{propos}
\noindent {\bf Proof:}\quad Suppose that $x$ and $y$ are given by
$[\xi]_0\in \Xi^r_0 X$ and $[\eta]_0\in \Xi^s_0 X$ respectively.
Set $\nabla x= \omega_i\tens [\xi_i]$ and
$\nabla x= \phi_i\tens [\eta_i]$,  for all $ [\xi_i]_0 \in \Xi^r_0 X$ 
and $[\eta_i]\in \Xi^s_0 X$,
then
\begin{eqnarray*}
\extd(\xi\wedge\eta) &=& \extd\xi\wedge\eta+(-1)^r\,\xi\wedge\extd\eta \cr
&=& \iota_* \omega_i\wedge\xi_i\wedge\eta
+(-1)^r\,\xi\wedge\iota_*\phi_i\wedge\eta_i\ .\quad\square
\end{eqnarray*}

\begin{propos} For all $x\in H^*(\Xi^*_0 X)$ and $\omega,\phi\in\Omega^*B$,
$$
(\id\wedge\sigma)(\sigma(x\tens \omega)\tens\phi) \,=\, \sigma(x\tens
(\omega\wedge\phi)).
$$
\end{propos}
\noindent {\bf Proof:}\quad Set $x=[\xi]_0$. We will use explicit summations
in this proof.
  We obtain
$$
(\id\wedge\sigma)(\sigma([\xi]_0\tens \omega)\tens\phi) = \sum_i
(\id\wedge\sigma)(\omega'_i \tens [\xi'_i]_0\tens\phi) 
= \sum_{ij} \omega'_i \wedge \phi'_{ij}\tens  [\xi''_{ij}]_0\ ,
$$
where
$$
\sum_i [\iota_*\omega'_i\wedge\xi'_i]_{|\omega|} = (-1)^{|\xi||\omega|}\,
[\xi\wedge\iota_*\omega]_{|\omega|}\ ,\quad
\sum_j [\iota_*\phi'_{ij}\wedge\xi''_{ij}]_{|\phi|} = (-1)^{|\xi||\phi|}\,
[\xi'_i\wedge\iota_*\phi]_{|\phi|}\ .
$$
 From this we obtain
\begin{eqnarray*}
\sum_{ij}[\iota_*(\omega'_i\wedge\phi'_{ij})\wedge 
\xi''_{ij}]_{|\omega|+|\phi|}
  &=& \sum_i  (-1)^{|\xi||\phi|}\, [\iota_*(\omega'_i)\wedge
  \xi'_i\wedge\iota_*\phi]_{|\omega|+|\phi|} \cr
  &=& (-1)^{|\xi|(|\phi|+|\omega|)}\, [\xi\wedge\iota_*\omega
  \wedge\iota_*\phi]_{|\omega|+|\phi|}\ .\quad\square
\end{eqnarray*}

\medskip
The reader will recall that in the construction of the spectral 
sequence the vector
spaces
  $\Omega^n B\tens_B H^m(\Xi_0^* X)$ appear. This is not such a simple thing
  as a tensor product differential complex, as the derivative involves
  a connection which does not map $H^*(\Xi_0^* X)$ to itself. It is therefore
  not surprising that the product structure has to be rather more 
complicated than
  the graded tensor product. In fact we have already given all the ingredients
  required for the product, it only remains to state them in a more 
coherent manner:

\begin{defin} \label{proddeffi}
Take  $(E^m,\nabla)\in {}_A\mathcal{F}$ for all $m\ge 0$,
and suppose that each $E^m$ is an $A$-bimodule. Give $e\in E^m$ the 
grade $|e|=m$.
A {\it product structure} on this
family consists of

(1)\quad $A$-bimodule
  maps $\sigma:E^m\tens_A \Omega^n A \to  \Omega^n A\tens_A E^m$,

(2)\quad a product $\wedge:E^m\tens_A E^{m'}\to E^{m+m'}$,

\noindent
which satisfy the following conditions, for all $e,f\in E^*$ and 
$\xi,\eta\in\Omega^*A$:

(a)\quad the product $(\xi\tens e)\wedge(\eta\tens f)=(-1)^{|e||\eta|}\,
\xi\wedge\sigma(e\tens\eta)\wedge f$ on $\Omega^*A\tens E^*$ is associative;

(b)\quad $(\id\wedge\sigma)(\nabla e\tens\xi) +\sigma(e\tens\extd\xi)
\,=\,[\extd\tens\id+(-1)^{|\xi|}(\id\wedge\nabla)]\sigma(e\tens\xi)$;

(c)\quad $\nabla(e\wedge f)=\nabla e\wedge f
+(\sigma\wedge\id)(e\tens\nabla f)$;

(d)\quad $(\id\wedge\sigma)(\sigma(e\tens \xi)\tens\eta) \,=\, \sigma(e\tens
(\xi\wedge\eta))$.

\end{defin}

\begin{propos} In Definition~\ref{proddeffi} the derivative 
$\nabla^{[*]}$ is a graded derivation
over the given product structure on $\Omega^*A\tens E^*$, i.e.\
\begin{eqnarray*}
\nabla^{[*]}((\xi\tens e)\wedge(\eta\tens f)) \,=\,
\nabla^{[*]}(\xi\tens e)\wedge(\eta\tens f)+(-1)^{|\xi|+|e|}
(\xi\tens e)\wedge\nabla^{[*]}(\eta\tens f)\ .
\end{eqnarray*}
Thus there is an induced product structure on the cohomology,
$$\wedge:H^n(A,E^m,\nabla)\tens H^{n'}(A,E^{m'},\nabla)\to
H^{n+n'}(A,E^{m+m'},\nabla).$$
\end{propos}
\noindent {\bf Proof:}\quad Begin with
\begin{eqnarray*}
(-1)^{|e||\eta|}\,\nabla^{[*]}((\xi\tens e)\wedge(\eta\tens f)) &=&
(\extd\tens\id+(-1)^{|\xi|+|\eta|}(\id\wedge\nabla))
(\xi\wedge\sigma(e\tens\eta)\wedge f) \cr
&=& \extd\xi\wedge\sigma(e\tens\eta)\wedge f
+(-1)^{|\xi|}\xi\wedge(\extd\tens\id)\sigma(e\tens\eta)\wedge f \cr
&&+\,(-1)^{|\xi|+|\eta|}\xi\wedge(\id\wedge\nabla)\sigma(e\tens\eta)\wedge 
f \cr
&&+\,(-1)^{|\xi|+|\eta|}\xi\wedge (\id\wedge\sigma\wedge\id)
(\sigma(e\tens \eta)\tens\nabla f)\ .
\end{eqnarray*}
Using property (b) of Definition~\ref{proddeffi}, this becomes
\begin{eqnarray} \label{kjud}
(-1)^{|e||\eta|}\,\nabla^{[*]}((\xi\tens e)\wedge(\eta\tens f)) &=& 
\extd\xi\wedge\sigma(e\tens\eta)\wedge f
+(-1)^{|\xi|}\xi\wedge(\id\wedge\sigma)(\nabla e\tens\eta)\wedge f \cr
&&+\,(-1)^{|\xi|}\xi\wedge\sigma(e\tens\extd\eta)\wedge f \cr
&&+\,(-1)^{|\xi|+|\eta|}\xi\wedge (\id\wedge\sigma\wedge\id)
(\sigma(e\tens \eta)\tens\nabla f)\ .
\end{eqnarray}
Next we calculate
\begin{eqnarray*}
(-1)^{|e||\eta|}\,\nabla^{[*]}(\xi\tens e)\wedge(\eta\tens f) &=&
(-1)^{|e||\eta|}\,(\extd\xi\tens e+(-1)^{|\xi|}\xi\wedge\nabla e)
\wedge(\eta\tens f)\ ,
\end{eqnarray*}
which is the same as the first two terms of (\ref{kjud}). Next
\begin{eqnarray*}
(-1)^{|e||\eta|+|\xi|+|e|}(\xi\tens 
e)\!\!\!\!\!\!\!\!&&\!\!\!\!\!\!\!\wedge\nabla^{[*]}(\eta\tens f)\cr
&=&
(-1)^{|e||\eta|+|\xi|+|e|}(\xi\tens e)\wedge(\extd \eta\tens f+(-1)^{|\eta|}
\eta\wedge\nabla f) \cr
&=& (-1)^{|\xi|} \xi\wedge\sigma(e\tens\extd\eta)\wedge f
+\, (-1)^{|\xi|+|\eta|} \xi\wedge 
(\sigma\wedge\id)(e\tens\eta\wedge\nabla f)\ ,
\end{eqnarray*}
so to prove the result we only need to verify
\begin{eqnarray*}
  (\id\wedge\sigma\wedge\id)
(\sigma(e\tens \eta)\tens\nabla f) &=&
  (\sigma\wedge\id)(e\tens\eta\wedge\nabla f)\ ,
\end{eqnarray*}
which is given by property (d) of Definition~\ref{proddeffi}.\quad$\square$

\section{Coactions of Hopf algebras}
In classical topology, fibrations arise whenever there is a 
continuous (compact) group action on a (compact) Hausdorff space 
(e.g.\ a free action gives rise to a principal fibration). A base of 
the fibration is then identified with the quotient of the total space 
by this action. In non-commutative geometry this corresponds to a 
coaction of a Hopf algebra on an algebra. This is the case that we 
consider in this section and, indeed in all the remaining sections.
\subsection{Differential calculi on Hopf algebras}
For more details on this subject, the reader should see \cite{worondiff}.
Suppose that a Hopf algebra $H$
with coproduct $\Delta_H$, counit $\epsilon_H$ and the invertible 
antipode $S$ has a differential calculus $\Omega^*H$.
We write the coproduct in $H$ as $\Delta_H(h) = h\sw 1\ot h\sw 2$, 
$\Delta_H^2(h) = h\sw 1\ot h\sw 2\ot h\sw 3$, etc., and the left 
$H$-coaction on $\Omega^*H$ as $\xi\mapsto
\xi_{[-1]}\tens \xi_{[0]}$ (summation understood). If there is no 
danger of confusion we will simply write $\Delta$ and  $\epsilon$ for 
$\Delta_H$ and $\epsilon_H$ (this convention applies to all other 
Hopf algebras as well).
In this section we shall not assume that
the coproduct is differentiable (this would give a bicovariant 
calculus), but only that
the left $H$-coaction $\lambda:\Omega^* H\to H\tens \Omega^* H$
is defined.
 $L^n H$ denotes the space of left invariant $n$-forms on $H$,
that is:
\begin{eqnarray*}
L^n H \,=\, {}^{{\rm co}H}(\Omega^n H)\,=\,\big\{ \xi\in \Omega^n H :
\xi_{[-1]}\tens \xi_{[0]}=1_H\tens \xi\big\}\ .
\end{eqnarray*}
The Hopf-Lie algebra
$\mathfrak{h}$ of $H$ is defined to be
\begin{eqnarray*}
\mathfrak{h} \,=\,\{\alpha:\Omega^1H\to k : 
\alpha(\eta\,h)=\alpha(\eta)\,\epsilon(h)
\quad\forall \eta\in\Omega^1 H\ ,\ \forall h\in H\}\ .
\end{eqnarray*}
Note that defining $\mathfrak{h}$ as a vector space only requires a 
`classical point',
that is an algebra map $\epsilon:H\to k$.

\begin{lemma} \label{dbvakldvbjh}
If, for a left invariant $\eta\in \Omega^1 H$, $\alpha(\eta)=0$ for 
all $\alpha\in\mathfrak{h}$,
then $\eta=0$.
\end{lemma}
\noindent {\bf Proof:}\quad For any $k$-linear map $T:L^1 H\to k$, define
$\alpha_T:\Omega^1H\to k$ by $\alpha_T(\xi) = 
T(\xi_{[0]}\,S^{-1}(\xi_{[-1]}))$. Then
for $h\in H$,
\begin{eqnarray*}
\alpha_T(\xi\,h) \,=\, 
T(\xi_{[0]}\,h_{(2)}\,S^{-1}(h_{(1)})\,S^{-1}(\xi_{[-1]}))\,=\,
\alpha_T(\xi)\,\epsilon(h)\ ,
\end{eqnarray*}
so $\alpha_T\in \mathfrak{h}$. For a left invariant $\eta\in \Omega^1 
H$, choose
$T$ so that $T(\eta)\neq 0$, and then $\alpha_T(\eta)\neq 0$. \quad$\square$

\subsection{Differentiable right coactions}

Suppose that the algebra $X$ has a differentiable right coaction $\rho$
(written on elements as  $\rho(x)=x_{[0]}\tens x_{[1]}\in X\tens H$, 
summation understood)
by the Hopf algebra $H$ which makes it into a comodule algebra.
This means that $\rho:X\to X\tens H$
is a coaction and a differentiable algebra map, so we obtain
a map of differential graded algebras (under the $\wedge$ multiplication)
\begin{eqnarray} \label{dirsum}
\rho_*:\Omega^n X \to \Omega^n(X\tens H) \,=\,
\bigoplus_{0\le r\le n} \Omega^r X \tens \Omega^{n-r}H\ .
\end{eqnarray}
Write $\Pi_{m,n-m}$ for the corresponding projection from 
$\Omega^n(X\tens H)$ to
$\Omega^m X\tens \Omega^{n-m}H$. Note that the maps $\Pi_{n,0}\rho_*$ 
define the right coactions of $H$ on $\Omega^nX$. These are also denoted 
by $\rho$. The subalgebra $B\subset X$ is defined to be the 
co-invariants for the right $H$-coaction, i.e.\ $B= X^{coH} := \{b\in 
X\; |\; \rho(b) = b\ot 1_H\}$.
We now define the calculus on $B$ by $\Omega^1 B=B.\extd B\subset 
\Omega^1 X$ and
$\Omega^n B=\bigwedge^n \Omega^1 B\subset \Omega^n X$. It is immediate that
$\Omega^n B \subset (\Omega^n X)^{{\rm co}H}$, the $H$-invariant 
$n$-forms on $X$. However we can be rather more restrictive:

\begin{defin} \label{ikyugc}
Define
\[
\mathcal{H}^n X=
\bigcap_{n> m\ge 0} \ker(
\Pi_{m,n-m}\rho_*:\Omega^n X\to \Omega^m X\tens \Omega^{n-m}H)\ .
\]
 The elements of $\mathcal{H}^n X$ are called  {\em horizontal} $n$-forms.
\end{defin}

\begin{remark}
It is immediate that $\Omega^n B \subset \mathcal{H}^n X$, and we might
conjecture that in `nice' cases we should have
$\Omega^n B = (\Omega^n X)^{{\rm co}H}\cap \mathcal{H}^n X$.
The reader should note that in the case of a bicovariant calculus
on $H$, the differential algebra $\Omega^* H$ is itself a graded Hopf algebra
(see \cite{Brz:rem}),
and then the conjecture is that $\Omega^* B$ is the invariant
part of $\Omega^* X$ under the right $\Omega^* H$-coaction.
\end{remark}

\begin{remark}
As in the classical case, it is possible to define horizontal 1-forms 
with reference
to the Hopf-Lie algebra.
Remember from \cite{braexp}
  that the vector fields on $X$
are the right $X$ module maps from $\Omega^1 X$ to $X$.
Every $\alpha\in\mathfrak{h}$ gives a vector field $\hat\alpha$ on 
$X$ defined by
$\hat\alpha(\xi)=(\id\tens\alpha)\Pi_{0,1}\rho_*(\xi)$, for every 
$\xi\in\Omega^1 X$.
\end{remark}

\begin{propos}
$\mathcal{H}^1 X=\cap_{\alpha\in\mathfrak{h}} 
\ker(\hat\alpha:\Omega^1 X\to X)$.
\end{propos}
\noindent {\bf Proof:}\quad First the reader should recall the 
definition of the cotensor product
 $U\, \square_H\, V$ of a right 
$H$-comodule $U$ and a left $H$-comodule $V$
 \cite{MilMoo}. This is 
the subset of $U\tens V$ consisting of all $u\tens v$ (summation 
implicit)
 where $u_{[0]}\tens u_{[1]}\tens v=u\tens v_{[-1]}\tens 
v_{[0]} \in U\tens H\tens V$. 
 Note that we can restrict the 
codomain to get $\Pi_{0,1}\rho_*:\Omega^1 X \to X \,
 \square_H\, 
\Omega^1 H$. Now there is a 1-1 correspondence between 
 $X \, 
\square_H\, \Omega^1 H$ and $X\tens L^1 H$ given by $x\, \square_H\, 
\xi \mapsto
 x_{[0]}\tens \xi\, S^{-1}(x_{[1]})$ and 
$y\tens\eta\mapsto y_{[0]} \, \square_H\, \eta\, y_{[1]}$. 
 This 
combines with Lemma \ref{dbvakldvbjh} to prove the 
result.\quad$\square$

\subsection{When the algebra coacted on is a Hopf algebra} 
\label{kausvc}
 A special case of interest, corresponding to 
homogenous spaces, is when the algebra $X$
 is itself a Hopf algebra.
Suppose that the Hopf algebra $X$ has a differentiable right coaction $\rho$
of the Hopf algebra $H$ which makes it into a comodule algebra.
We shall also assume that  $\rho$
commutes with the coproduct $\Delta_X$ of $X$, i.e.\ 
$(\id\tens\rho)\Delta_X=(\Delta_X\tens\id)\rho:
X\to X\tens X\tens H$. This is the case if and only if the map
$\pi:X\to H$ defined by $\pi(x)=(\epsilon_X\tens\id)\rho(x)$ is a bialgebra
map, and then $\rho(x)=(\id\tens\pi)\Delta_X(x)$.

\begin{defin}[\cite{Swe:int}]
A {\it normalised left integral} for $H$ is a map
$\int:H\to k$ with $(\id\tens\int)\Delta_H=1_H.\int:H\to H$
and $\int 1_H=1$.
\end{defin}

One easily checks that if  $\int$ is a normalised left integral for $H$, then the map $(\id\tens\int)\rho:X\to X$
is a projection onto the co-invariant subalgebra $B= X^{coH}$, for
any right $H$-comodule algebra $X$. Furthermore, if $X$ is itself a Hopf algebra
and $\rho$ commutes with $\Delta_X$, then we can state the following

\begin{lemma}\label{oiubdwijkuvb} If $H$ has a normalised left integral, then
$\Delta_X B\subset X\tens B$.
\end{lemma}
\noindent {\bf Proof:}\quad By the definition of co-invariants, for 
all $b\in B$,  $\rho(b)=b\tens 1_H$,
so
\begin{eqnarray*}
(\Delta_X\tens\id)\rho(b)\,=\,b_{(1)}\tens b_{(2)}\tens 
1_H\,=\,b_{(1)}\tens\rho(b_{(2)})\ .
\end{eqnarray*}
If we apply the integral to this we obtain
\begin{eqnarray*}
b_{(1)}\tens b_{(2)}\,=\,b_{(1)}\tens(\id\tens\int)\rho(b_{(2)}) \in 
X\tens B\ .\quad\square
\end{eqnarray*}
 
 \medskip This means that the Hopf algebra $X$ left 
coacts on $B$. Thus $B$ can be viewed as a noncommutative 
generalisation of a homogenous space of $X$ \cite{Pod:sph}.

\section{The non-commutative Hopf fibration with a non-bi\-co\-variant calculus}
In this section we give an explicit example of a non-commutative 
differentiable fibration. It is well known that the underlying 
algebra inclusion is a {\em quantum principal bundle} 
\cite{BrzMaj:gau}, our aim, however, is to show that it is a
differentiable fibration in the sense of Definition~\ref{jasdfvajkhv}.

\subsection{Example: The quantum Hopf fibration} \label{kuasceewmnoon}
This is an example of the type of coaction discussed in Section~\ref{kausvc}.
Consider the complex Hopf algebra $X=\mathcal{A}(SL_q(2))$ generated by
$\{\alpha,\beta,\gamma,\delta\}$ with the relations
\begin{eqnarray}
&& \alpha\beta\,=\, q\,\beta\alpha\ ,\quad
\alpha\gamma\,=\, q\,\gamma\alpha\ ,\quad
\beta\gamma\,=\, \gamma\beta\ ,\quad
\beta\delta\,=\, q\,\delta\beta\ ,\quad
\gamma\delta\,=\, q\,\delta\gamma\ , \cr
&&\alpha\delta\,=\, \delta\alpha+(q-q^{-1})\,\beta\gamma\ ,\quad
\alpha\delta-q\,\beta\gamma\,=\,1\ ,
\end{eqnarray}
where $q$ is a complex number which is not a root of unity. On this 
level of algebraic generality, there is no need to make further 
restrictions on $q$, although geometrically
most interesting is the case $0<q<1$, whereby $X$ can be made into a 
$*$-algebra and extended to a $C^*$-algebra of functions on the 
quantum group $SU_q(2)$
(cf.\ \cite{Wor:twi}).
The coproduct is given by
\begin{eqnarray}
&&\Delta\alpha\,=\,\alpha\tens\alpha+\beta\tens\gamma\ ,\quad
\Delta\beta\,=\,\alpha\tens\beta+\beta\tens\delta\ ,\cr
&& \Delta\gamma\,=\,\gamma\tens\alpha+\delta\tens\gamma\ ,\quad
\Delta\delta\,=\, \delta\tens\delta+\gamma\tens\beta\ ,
\end{eqnarray}
and counit and antipode by
\begin{eqnarray*}
&& \epsilon(\alpha)\,=\,\epsilon(\delta)\,=\,1\ ,\quad
\epsilon(\beta)\,=\,\epsilon(\gamma)\,=\,0\ ,\cr
&& S(\alpha)\,=\,\delta\ ,\quad
S(\delta)\,=\,\alpha\ ,\quad
S(\beta)\,=\,-q^{-1}\,\beta\ ,\quad
S(\gamma)\,=\,-q\,\gamma\ .
\end{eqnarray*}

We will take $H$ to be the group algebra of $\mathbb{Z}$, which we 
take as generated by
$z$, $z^{-1}$ with $\Delta z^{\pm 1}=z^{\pm 1}\tens z^{\pm 1}$,
$S(z^{\pm 1})=z^{\mp 1}$ and $\epsilon(z^{\pm 1})=1$.
The Hopf algebra map $\pi:X\to H$ is given by
\begin{eqnarray*}
\pi(\alpha)\,=\,z\ ,\quad
\pi(\delta)\,=\,z^{-1}\ ,\quad
\pi(\beta)\,=\,\pi(\gamma)\,=\,0\ .
\end{eqnarray*}
The right $H$-coaction $\rho$ on $X$ is then given by
\begin{eqnarray*}
\rho(\alpha)\,=\,\alpha\tens z\ ,\quad
\rho(\beta)\,=\,\beta\tens z^{-1}\ ,\quad
\rho(\gamma)\,=\,\gamma\tens z\ ,\quad
\rho(\delta)\,=\,\delta\tens z^{-1}\ .
\end{eqnarray*}
The invariant part of $X$, $B=X^{{\rm co}H}=\mathcal{A}(S^2_q)$ is 
generated as an algebra by
$\{\alpha\beta,\alpha\delta,\gamma\delta\}$ and is known as (an algebra of
functions on) the  standard quantum 2-sphere \cite{Pod:sph}.

\subsection{The 3D non-bicovariant calculus on $\mathcal{A}(SL_q(2))$}
This left covariant differential calculus on $X=\mathcal{A}(SL_q(2))$ 
was introduced by
Woronowicz in \cite{Wor:twi} and is generated by 3 left invariant 1-forms
$\{\omega^0,\omega^1,\omega^2\}$. The differentials of the generators 
are given by
\begin{eqnarray} \label{3Drel1}
\extd \alpha \,=\, \alpha\,\omega^1-q\,\beta\,\omega^2\ ,\quad
\extd \beta \,=\, \alpha\,\omega^0-q^2\,\beta\,\omega^1\ ,\cr
\extd \gamma \,=\, \gamma\,\omega^1-q\,\delta\,\omega^2\ ,\quad
\extd \delta \,=\, \gamma\,\omega^0-q^2\,\delta\,\omega^1\ .
\end{eqnarray}
We have the commutation relations
\begin{eqnarray}\label{3Drel8}
\omega^0\,\alpha \,=\, q^{-1}\,\alpha\,\omega^0 &,&
\omega^0\,\beta \,=\, q\,\beta\,\omega^0\ ,\cr
\omega^1\,\alpha \,=\, q^{-2}\,\alpha\,\omega^1 &,&
\omega^1\,\beta \,=\, q^2\,\beta\,\omega^1\ ,\cr
\omega^2\,\alpha \,=\, q^{-1}\,\alpha\,\omega^2 &,&
\omega^2\,\beta \,=\, q\,\beta\,\omega^2\ ,
\end{eqnarray}
and similarly for replacing $\alpha\to\gamma$ and $\beta\to\delta$.
For the higher forms we have exterior derivative
\begin{eqnarray}\label{3Drel9}
\extd \omega^0 \,=\,q^2(q^2+1)\,\omega^0\wedge\omega^1\ ,\quad
\extd \omega^1 \,=\,q\,\omega^0\wedge\omega^2\ ,\quad
\extd \omega^2 \,=\,q^2(q^2+1)\,\omega^1\wedge\omega^2\ ,
\end{eqnarray}
and wedge multiplication
\begin{eqnarray}\label{3Drel10}
&&\omega^0\wedge\omega^0\,=\,\omega^1\wedge\omega^1\,=\,
\omega^2\wedge\omega^2 \,=\, 0\ ,\cr
&& \omega^2\wedge\omega^0\,=\,-q^2\,\omega^0\wedge\omega^2\ ,\quad
\omega^1\wedge\omega^0\,=\,-q^4\,\omega^0\wedge\omega^1\ ,\quad
\omega^2\wedge\omega^1\,=\,-q^4\,\omega^1\wedge\omega^2\ .
\end{eqnarray}

\subsection{The differentiable coaction} \label{yuasdfg}
We need the map $\pi$ given in
Section~\ref{kuasceewmnoon} to extend to a map $\pi_*$ of 
differential graded algebras. Such an extension of $\pi$ exists, 
provided there is a suitable differential structure
on $H$, which can be constructed as follows.
{}From (\ref{3Drel1}) we obtain $\extd z=z\,\pi_*(\omega^1)$, 
$0=z\,\pi_*(\omega^0)$,
$0=-q\,z^{-1}\,\pi_*(\omega^2)$ and 
$\extd(z^{-1})=-q^2\,z^{-1}\,\pi_*(\omega^1)$.
This can be summarised by
\begin{eqnarray}
\pi_*(\omega^0)=\pi_*(\omega^2)= \, 0\ ,\quad \pi_*(\omega^1)=z^{-1}.\extd z\ ,
\quad z.\extd z = q^2\,\extd z.z\ .
\end{eqnarray}
(To see this, note that from $z.z^{-1}=1$ we use the derivation 
property for $\extd$
to get $\extd (z^{-1})=-z^{-1}.\extd z.z^{-1}$.)
It is easily checked that the map $\pi_*$ defined in this fashion 
satisfies all the relations and that the constructed differential 
calculus on $H$ is bicovariant. However the cost of differentiability 
of $\pi_*$ is that the commutative algebra
$H$ is given a noncommutative differential structure!

To find $\rho_*$ we look at (\ref{3Drel1}), and use $\rho_*(\extd\alpha)
=\extd(\rho(\alpha))$ etc.\ to give
\begin{eqnarray} \label{3Drel11}
\rho_*(\omega^0)\,=\, \omega^0\tens z^{-2}\ ,\quad
\rho_*(\omega^1)\,=\, 1\tens z^{-1}.\extd z +\omega^1\tens 1\ ,\quad
\rho_*(\omega^2)\,=\, \omega^2\tens z^{2}\ .
\end{eqnarray}
To check that this gives a well defined map on $\Omega^1 X$, one needs to
check that it is consistent with the relations in (\ref{3Drel8}) -- 
this is left to the reader. Then to define
$\rho_*$ on the higher forms by using the wedge product we only have to check
the relations in (\ref{3Drel9}) and (\ref{3Drel10}), which is easily 
done by a straightforward calculation.

To find the horizontal 1-forms we apply $\Pi_{0,1}$ to (\ref{3Drel11}) to get
\begin{eqnarray*}
\Pi_{0,1}\rho_*(\omega^1)\,=\,1\tens z^{-1}\extd z\ ,\quad
\Pi_{0,1}\rho_*(\omega^0)
\,=\,\Pi_{0,1}\rho_*(\omega^2)\,=\,0\  .
\end{eqnarray*}
It follows that the horizontal 1-forms are precisely those of the form
$a\,\omega^0+b\,\omega^2$ for $a,b\in X$.
We can also calculate the right $H$-coaction
by applying $\Pi_{1,0}$ to (\ref{3Drel11}) to get
\begin{eqnarray*}
\Pi_{1,0}\rho_*( \omega^1)\,=\,\omega^1\tens 1\ ,\quad
\Pi_{1,0}\rho_* (\omega^0)\,=\,\omega^0\tens z^{-2}\ ,\quad
\Pi_{1,0}\rho_* (\omega^2)\,=\,\omega^2\tens z^2\ .
\end{eqnarray*}
Then the invariant horizontal 1-forms are precisely those of the form
$a\,\omega^0+b\,\omega^2$ where $\rho(a)=a\tens z^2$ and
$\rho(b)=b\tens z^{-2}$.

\subsection{The corresponding calculus on $B=\mathcal{A}(S^2_q)$}
We can calculate
\begin{eqnarray}\label{skjcvg}
\extd(\alpha\beta) &=& \alpha^2\,\omega^0 - q^2\,\beta^2\,\omega^2\ ,\cr
q\,\extd(\beta\gamma) &=& \alpha\gamma\,\omega^0 - 
q^2\,\beta\delta\,\omega^2\ ,\cr
\extd(\gamma\delta) &=& \gamma^2\,\omega^0 - q^2\,\delta^2\,\omega^2\ .
\end{eqnarray}
 From this we get
\begin{eqnarray*}
\delta\,\extd(\alpha\beta) - q^{-1}\,\beta\,\extd(\beta\gamma) &=& 
\alpha\,\omega^0\ ,\cr
q\,\delta\,\extd(\beta\gamma)-q^{-1}\,\beta\,\extd(\gamma\delta) &=& 
\gamma\,\omega^0\ .
\end{eqnarray*}
By left multiplying these last equations by $\alpha$ and $\gamma$ we see
that $\alpha^2\,\omega^0$, $\alpha\gamma\,\omega^0$ and $\gamma^2\,\omega^0$
are all in $B.\extd B$. From (\ref{skjcvg}) we deduce that
$\beta^2\,\omega^2$, $\beta\delta\,\omega^2$ and $\delta^2\,\omega^2$
are also all in $B.\extd B$.

Given a monomial $a$ in the generators $\{\alpha,\beta,\gamma,\delta\}$
with $\rho(a)=a\tens z^2$, we can reorder it as either
$a=x\,\alpha^2$ or $a=x\,\alpha\gamma$ or $a=x\,\gamma^2$, where
$x\in B$. Thus we have $a\,\omega^0\in B.\extd B$.  Likewise for a monomial
$b$ with $\rho(b)=b\tens z^{-2}$ we have $b\,\omega^2\in B.\extd B$.
 From this and the discussion in Section~\ref{yuasdfg}
  we conclude that $\Omega^1 B$ is precisely the horizontal
invariant 1-forms on $X$.

Now we shall consider the 2-forms. Since $\rho_*$ is a graded algebra 
map, we immediately obtain
$$
\rho_*(\omega^0\wedge\omega^2) = \omega^0\wedge\omega^2\otimes 1,
\qquad \rho_*(\omega^l\wedge\omega^1) = \omega^l \otimes 
z^{2l-3}.\extd z + \omega^l\wedge\omega^1\otimes z^{2l-2}, \qquad l 
=0,2.
$$
  Hence
the horizontal 2-forms are multiples of $\omega^0\wedge\omega^2$.
Then the invariant horizontal 2-forms are $B.\omega^0\wedge\omega^2$.
To see that $\Omega^2 B$ is all of this, we use the following relation:
\begin{eqnarray*}
\alpha^2\delta^2 - 
(q+q^{-1})\,\alpha\gamma\beta\delta+q^2\,\gamma^2\beta^2\,=\,1\ .
\end{eqnarray*}
By using $\alpha^2\,\omega^0\wedge\delta^2\,\omega^2=q^2\,\alpha^2
\delta^2\,\omega^0\wedge\omega^2$ and similar calculations,
we see that $\omega^0\wedge\omega^2$ is contained in
$\Omega^1 B\wedge \Omega^1 B$.

All 3-forms are multiples of $\omega^0\wedge\omega^1\wedge\omega^2$,
but none of these (except zero) are horizontal, so we conclude that
$\Omega^3 B=0$.

\subsection{An easy example of a spectral sequence}
We will use the notation $\<\dots\>$ to denote the right $X$-module
generated by the listed elements. Then as right $X$-modules, 
$B\tens_B X\cong X$,
$\Omega^1 B\tens_B X\cong \<\omega^0,\omega^2\>$ and
$\Omega^2 B\tens_B X\cong \<\omega^0\wedge\omega^2\>$. We can calculate
the $\Xi_m^n X$ as shown in the following table:

\begin{center}
{\begin{tabular}{c|c|c|c}$\Xi^n_m X$ & $n=0$ & $n=1$ & $n>1$ \\\hline 
$m=0$ & $X$ & $\<\omega^1\>$ & 0 \\\hline $m=1$ & 
$\<\omega^0,\omega^2\>$ & 
$\<\omega^0\wedge\omega^1,\omega^2\wedge\omega^1\>$ & 0 \\\hline 
$m=2$ & $\<\omega^0\wedge\omega^2\>$ & 
$\<\omega^0\wedge\omega^1\wedge\omega^2\>$ & 0 \\\hline $m>2$ & 0 & 0 
& 0\end{tabular}}
\end{center}
\smallskip

\noindent It follows that
$\Theta_m:\Omega^m B \tens_B\Xi_0^n X=\Omega^m B \tens_B X\tens_X \Xi_0^n X
\to \Xi_m^n X$ (as defined in Definition~\ref{hhhddd}) is an isomorphism, and that
  the quantum Hopf fibration $\iota: \mathcal{A}(S_q^2)\hookrightarrow 
\mathcal{A}(SL_q(2))$ is a differential fibration for this 
differential structure.
  \smallskip

\noindent Now we shall calculate the $E_2$ page of the spectral 
sequence in this case.
The first thing to do is to look at $H^*(\Xi_0^* X)$. Recall that we 
consider only the generic case, where
$q$ is not a root of unity. Note that the coaction $\rho$ makes $X$ a 
$\Z$-graded
algebra with the grading $\deg \alpha = \deg \gamma =1$, $\deg\beta = 
\deg\delta =-1$, $\deg 1 = 0$.

\begin{lemma} \label{ksahn}
For any homogeneous $x\in X$, the differential $\extd:\Xi_0^0 X=X \to 
\Xi_0^1 X=\Omega^1 X/\<\omega^0,\omega^2\>$
gives
$$
\extd x \,=\, [\deg x;q^{-2}]\,x\,\omega^1,
$$
where $[n;q^{-2}] = \frac{q^{-2n}-1}{q^{-2}-1}$ is a $q^{-2}$-integer.
\end{lemma}
\noindent {\bf Proof:}\quad This is most easily proved by using the 
linear basis
 $\{\alpha^a\beta^b\gamma^c,\beta^b\gamma^c\delta^d\}$
 of $\mathcal{A}(S^2_q)$. 
\quad$\square$

\begin{propos}
As left $B$-modules, $H^0(\Xi_0^* X)=B$, $H^1(\Xi_0^* X)=B.\omega^1$,
and for $n>1$, $H^n(\Xi_0^* X)=0$.
\end{propos}
\noindent {\bf Proof:}\quad This comes from Lemma \ref{ksahn}
and $\Xi_0^n X=0$ for $n>1$. \quad$\square$

\begin{remark}
We now have to find the left $B$-connection $\nabla$ described in 
Proposition \ref{uioy}.
As each $H^n(\Xi_0^* X)$ is a finitely generated $B$-module,
it is enough to find $\nabla$ on the generators. Choose
generators $1_B$ and $\omega^1$ in $H^0(\Xi_0^* X)$ and $H^1(\Xi_0^* X)$
respectively; an explicit calculation then implies that $\nabla 1_B=0$
and $\nabla \omega^1=0$. Now we can calculate the $\nabla$-cohomology
of the $H^n(\Xi_0^* X)$ module, which is given by the cochain complex
\begin{eqnarray*}
H^n(\Xi_0^* X) \to \Omega^1 B\tens_B H^n(\Xi_0^* X)
\to \Omega^2 B\tens_B H^n(\Xi_0^* X) \to\dots\ .
\end{eqnarray*}
Using the generators, we identify this with the usual de Rham complex
\begin{eqnarray*}
B \to \Omega^1 B
\to \Omega^2 B \to\dots\ ,
\end{eqnarray*}
and so we get $E_2^{p,r}\cong H^p_{dR}(B)$ for $r=0,1$, and
$E_2^{p,r}\cong 0$ for other values of $r$. This gives the $E_2$ page 
of the Serre
spectral sequence (we display only potentially non-zero terms)
\newline\noindent
\begin{picture}(200,100)(-60,-20)
\put(0,0){\vector(1,0){270}}
\put(0,0){\vector(0,1){70}}
\put(260,-9){$p$}
\put(-12,60){$r$}

\put(-10,10){$0$}
\put(-10,40){$1$}

\put(10,10){$H^0_{dR}(B)$}
\put(10,40){$H^0_{dR}(B)$}
\put(21,-10){$0$}
\put(60,10){$H^1_{dR}(B)$}
\put(60,40){$H^1_{dR}(B)$}
\put(71,-10){$1$}
\put(110,10){$H^2_{dR}(B)$}
\put(110,40){$H^2_{dR}(B)$}
\put(121,-10){$2$}
\put(160,10){$H^3_{dR}(B)$}
\put(160,40){$H^3_{dR}(B)$}
\put(171,-10){$3$}
\end{picture}
\newline\noindent
The only possibly non-zero differentials on this page are
$\extd_2:(0,1)\to(2,0)$ and $\extd_2:(1,1)\to(3,0)$. All further pages
have differentials all zero, just from considering the
indices. From this we see that $H^3_{dR}(B)\cong H^4_{dR}(X)$, but
$H^4_{dR}(X)=0$ as $\Omega^4 X=0$, so $H^3_{dR}(B)=0$.
Using this, we get $H^2_{dR}(B)\cong H^3_{dR}(X)$. Also we obtain
$H^0_{dR}(B)\cong H^0_{dR}(X)$ and the more complicated
cases
\begin{eqnarray*}
H^1_{dR}(X) &\cong& H^1_{dR}(B) \oplus \ker(\extd_2:
H^0_{dR}(B)\to H^2_{dR}(B))\ , \cr
H^2_{dR}(X) &\cong& H^1_{dR}(B) \oplus {\rm coker}(\extd_2:
H^0_{dR}(B)\to H^2_{dR}(B))\ .
\end{eqnarray*}
To get any further, we would have to use additional information about
either $B$ or $X$. However this is one of the primary reasons why the 
Serre spectral sequence is useful, it turns information about one 
space into information about the other space.

\end{remark}

\section{A construction for bicovariant calculi}\label{sec.bicov}
In this section we consider Hopf algebras
$X$ and $H$ with bicovariant differential calculi. We assume  that there
exists a differentiable surjective Hopf algebra map $\pi:X\to H$. The 
right $H$-coaction
on $X$ is given by $\rho=(\id\tens\pi)\Delta:X\to X\tens H$ (cf.\ 
Section~\ref{kausvc}).
Since the calculus on $X$ is bicovariant the coproduct $\Delta$ in $X$ is a
differentiable map, hence also the coaction $\rho$ is differentiable (as a
composition of differentiable maps).

\subsection{Left invariant forms and coactions}
We first study the 
covariance properties of the spaces of horizontal $n$-forms
(see 
Definition \ref{ikyugc}).

\begin{propos} \label{isoufxgh}
$\mathcal{H}^n X$ is preserved by the right $H$-coaction,
i.e.\ $\rho(\mathcal{H}^n X)\subset \mathcal{H}^n X \tens H$.
\end{propos}
\noindent {\bf Proof:}\quad Start with any $\eta\in \mathcal{H}^n X$.
To check that $\rho(\eta)\in \mathcal{H}^n X \tens H$ we
need to show that $(\Pi_{m,n-m}\rho_*\tens\id)\rho(\eta)=0$,
for all $n> m\ge 0$. Inserting the definition of the right
coaction, we need to show that 
$(\Pi_{m,n-m}\rho_*\tens\id)\Pi_{n,0}\rho_*(\eta)=0$
for all $n> m\ge 0$. This, using more projections to forms,
is the same as $\Pi_{m,n-m,0}(\rho_*\tens \id)\rho_*(\eta)=0$
(here 
we have extended the projection $\Pi$ to three factors in the obvious 
manner). By
the coaction property this is $\Pi_{m,n-m,0}(\id\tens\Delta_*)\rho_*(\eta)=0$.
However as $\eta\in \mathcal{H}^n X$ we know that
$\rho_*(\eta)\in \Omega^n X\tens H$, giving
$(\id\tens\Delta_*)\rho_*(\eta)\in\Omega^n X\tens H\tens H$, and 
applying the projection gives zero.\quad$\square$

\begin{propos} \label{leftisoufxgh}
The space of horizontal $n$-forms $\mathcal{H}^n X$ is preserved by
the left $X$-coaction,
i.e.\ $\Pi_{0,n}\Delta_*(\mathcal{H}^n X)\subset X\tens \mathcal{H}^n X$.
\end{propos}
\noindent {\bf Proof:}\quad For $\eta\in \mathcal{H}^n X$ and all
$0\le m<n$ we need to show that
\begin{eqnarray*}
(\id\tens\id\tens\pi_*)(\id\tens\Pi_{m,n-m}\Delta_*)\Pi_{0,n}\Delta_*\eta
\end{eqnarray*}
vanishes. By coassociativity, this is the same as
\begin{eqnarray*}
\Pi_{0,m,n-m}(\Delta_*\tens\id)(\id\tens\pi_*)\Delta_*\eta\ .
\end{eqnarray*}
Now $(\id\tens\pi_*)\Delta_*\eta\in \Omega^n X\tens H$, so applying 
the projection
gives zero, as $n-m>0$.\quad$\square$

\subsection{The 1-forms on the base $B$}
To identify the differential 
forms on the base $B$, we require that $\pi:X\to H$
satisfies an 
additional condition, and this is best phrased in terms 
of the space
\begin{eqnarray}
\mathcal{K} \,=\, \ker(\pi_*:L^1 X\to L^1 H)\ 
.
\end{eqnarray}
Note that $\mathcal{K} $ is simply the space of horizontal left 
invariant 1-forms on $X$.

\begin{defin} \label{fcond}
We say that $\pi:X\to H$ {\em satisfies condition K} if $\mathcal{K} 
\subset \extd B.X$.
\end{defin}
Note that checking that $\pi$ satisfies condition K is easier than it 
might seem, as often the left invariant 1-forms on $X$ form a finite 
dimensional space (see
the explicit example in Section~\ref{sub.check.F}).

\begin{propos} \label{kszcvvbxcasjh}
If $\pi:X\to H$ satisfies condition K, then $\mathcal{H}^1 X=\extd B.X$.
\end{propos}
\noindent {\bf Proof:}\quad Take any $\eta\in\mathcal{H}^1 X$. By 
Proposition~\ref{leftisoufxgh} on the left $X$-coaction, 
$
\eta_{[-1]} \tens \eta_{[0][-1]} \tens\eta_{[0][0]} \in X\tens 
X\tens \mathcal{H}^1 X.
$
Remember for any 1-form $\xi$, that $\xi_{[0]}.S^{-1}(\xi_{[-1]})$ is left
invariant.
It follows that $\eta_{[-1]}  \tens\eta_{[0][0]}.S^{-1}(\eta_{[0][-1]}) \in
X\tens \extd B.X$, so
$$\eta=\eta_{[0][0]}.S^{-1}(\eta_{[0][-1]}) \eta_{[-1]}
\in \extd B.X\ .
$$
The other inclusion is immediate.\quad $\square$

\begin{propos}
Suppose that $H$ has a normalised left integral, and that $\pi:X\to 
H$ satisfies condition K.
Then $\Omega^1 B = \extd B.B=(\mathcal{H}^1 X)^{{\rm co}H}$.
\end{propos}
\noindent {\bf Proof:}\quad Apply the left integral to
the result of Proposition~\ref{kszcvvbxcasjh}.\quad$\square$

\subsection{Bicovariant calculi on Hopf algebras using left invariant 1-forms}
\label{isuydfgc}
In the case where the coproduct is differentiable, there is a 
construction of the
calculus on a Hopf algebra $X$ in terms of the left invariant 1-forms $L^1 X$
which is due to Woronowicz \cite{worondiff}.

There is an isomorphism of $X$ modules and comodules (the 
module/comodule structures indicated by the dots)
\begin{equation}\label{isom.1}
{}^\bullet\Omega^1 X {}_\bullet^\bullet \to
  (L^1 X)^\bullet \tens {}^\bullet X{}_\bullet^\bullet,
\qquad \xi \mapsto \xi_{[0]}\,S^{-1}(\xi_{[-1]}\sw 2)\tens \xi_{[-1]}\sw 1,
\end{equation}
  with the inverse given
by the product map. It is also a left $X$-module map, but with the left action
on $ L^1 X \tens X{}$ given by $x\, \la (\eta\tens y)=x_{(2)}\la\eta\tens 
x_{(1)}\,y$,
and $x\, \la\eta=x_{(2)}\,\eta\,S^{-1}(x_{(1)})$ for $\eta\in L^1 X$ and 
$x,y\in X$.

The relation between the left $X$-action on $L^1 X$ and the right $X$-coaction
$\rho_X: L^1 X\to L^1 X\tens X$ is summarised in
the equation $\rho_X(x\, \la\eta)=x_{(2)}\la \eta_{[0]}\tens
x_{(3)}\, \eta_{[1]}\,S^{-1}(x_{(1)})$. This fits the left action - 
right coaction
version of a Yetter-Drinfeld module (cf.\  \cite[Section~5]{CaMiZhFrobenius}).
By the standard results on Yetter-Drinfeld modules there is a braiding
$\sigma:L^1 X\tens L^1 X\to L^1 X\tens L^1X$ defined by
$\sigma(\xi\tens\eta)=\eta_{[0]}\tens S(\eta_{[1]})\la \xi$, with inverse
$\sigma^{-1}(\xi\tens\eta )= \xi_{[1]}\la \eta\tens \xi_{[0]}$.

We define the wedge product on $L^1X$ as a quotient
\begin{eqnarray*}
L^1X \wedge L^1X \,=\,\frac{L^1X \tens L^1X }
{\ker(\sigma-\id\tens\id:L^1X \tens L^1X\to L^1X \tens L^1X)}\ ,
\end{eqnarray*}
and extend this to higher wedge products.
There is (as a matter of definition
of the higher forms) an isomorphism
\begin{equation}\label{isom.2}
\Omega^n X \to (L^1 X)^{\wedge n}\tens X,
\end{equation}
with wedge product defined by $(\xi\tens x)\wedge(\eta\tens y)=
\xi\wedge(x_{(2)}\la\eta)\tens x_{(1)}\,y$.

\subsection{Identifying the higher dimensional calculus on the base 
algebra}

\begin{propos} \label{soidudvg}
Using the isomorphism \eqref{isom.1}, $\mathcal{H}^1 X$ corresponds to
$\mathcal{K}\tens X$.
\end{propos}
\noindent {\bf Proof:}\quad Begin with a horizontal one-form $\xi\in 
\mathcal{H}^1 X$, and apply the
isomorphism \eqref{isom.1} to get
$\xi_{[0]}\,S^{-1}(\xi_{[-1]}\sw 2)\tens \xi_{[-1]}\sw 1$. We need to show that
$\pi_*(\xi_{[0]})\,\pi(S^{-1}(\xi_{[-1]}\sw 2))\tens \xi_{[-1]}\sw 1=0$.
As $\xi\in \mathcal{H}^1 X$ we know that $\xi_{[-1]}\tens\pi_*(\xi_{[0]})=0$,
and the required result follows from this.

For the other direction, take $\eta\tens x\in \mathcal{K}\tens X$. Then
applying the left $X$-action to $\eta\,x$ gives $x_{(1)}\tens 
\eta\,x_{(2)}$, and applying $\id\tens\pi_*$ to this gives
$x_{(1)}\tens \pi_*(\eta)\,\pi(x_{(2)})=0$.
\quad $\square$

\begin{propos} \label{sjkhcv}
The usual right $H$-coaction on $L^1 X$ restricts to one on $\mathcal{K}$. Also
the usual left $X$-action on $L^1 X$ restricts to one on $\mathcal{K}$.
\end{propos}
\noindent {\bf Proof:}\quad For the coaction, for $\eta\in 
\mathcal{K}$ we need to show
that $\Pi_{1,0} (\id\tens\pi_*)\Delta_*\eta\in \mathcal{K}\tens H$. To do this
we need to show the vanishing of $\Pi_{1,0} (\pi_*\tens\pi_*)\Delta_*\eta
=\Pi_{1,0} \Delta_*\pi_*(\eta)=0$ (using the fact that $\pi$ is a 
coalgebra map).

For the action, we have
\begin{eqnarray*}
\pi_*(x\la\eta)\,=\,\pi_*(x_{(2)}\,\eta\,S^{-1}(x_{(1)}))\,=\,
\pi(x_{(2)})\,0\,\pi(S^{-1}(x_{(1)}))\,=\,0\ .\quad\square
\end{eqnarray*}

\begin{cor} \label{askfjascv}
If $\pi:X\to H$ satisfies condition K and $H$ has a normalised left integral,
  then, using the isomorphism \eqref{isom.1},
  $\Omega^1 B$ corresponds to $(\mathcal{K}\tens X)^{{\rm co}H}$.
\end{cor}

\begin{lemma} \label{askfjauuu}
If $H$ has a normalised left integral $\int$,
then for all $a,c\in H$
\begin{eqnarray*}
c_{(2)}\int\Big(a\,S(c_{(1)})\Big) \,=\, a_{(1)}\int\Big(a_{(2)}\,S(c)\Big)\ .
\end{eqnarray*}
\end{lemma}
\noindent {\bf Proof:}\quad For all $a,b\in H$ the left integral property gives
\begin{eqnarray*}
a_{(1)}\tens a_{(2)}b_{(1)}\int\Big(a_{(3)}b_{(2)}\Big) \,=\,
a_{(1)}\tens 1_H\int\Big(a_{(2)}b\Big)\ .
\end{eqnarray*}
Applying $S^{-1}$ to the last factor gives
\begin{eqnarray*}
a_{(1)}\tens S^{-1}(b_{(1)}) S^{-1}(a_{(2)})\int\Big(a_{(3)}b_{(2)}\Big) \,=\,
a_{(1)}\tens 1_H\int\Big(a_{(2)}b\Big)\ .
\end{eqnarray*}
Now multiply the second factor on the right by the first to get
\begin{eqnarray*}
  S^{-1}(b_{(1)}) \int\Big(a\,b_{(2)}\Big) \,=\,
a_{(1)}\int\Big(a_{(2)}b\Big)\ .
\end{eqnarray*}
Finally, putting $b=S(c)$ gives the result.\quad $\square$

\begin{lemma} \label{askfjascvuuu}
If $H$ has a normalised left integral,
  then there is a projection $p:L^1 X\to L^1 X$
  with image $\mathcal{K}$ which preserves the right $H$-coaction
  (i.e.\ $p$ is right $H$-colinear, i.e.\ $\rho\, p=(p\tens\id)\rho$).
\end{lemma}
\noindent {\bf Proof:}\quad Take any linear projection $p_0:L^1 X\to L^1 X$
  with image $\mathcal{K}$, and define (using square brackets for the 
$H$-coaction)
  \begin{eqnarray*}
p(\xi)\,=\,p_0(\xi_{[0]})_{[0]}\, 
\int\Big(p_0(\xi_{[0]})_{[1]}\,S(\xi_{[1]})\Big)\ .
\end{eqnarray*}
First we show that $p$ is a projection to $\mathcal{K}$. Since the image
of $p_0$ is $\mathcal{K}$, and $\mathcal{K}$ is coacted on by $H$, it is
obvious from the formula that the image of $p$ is contained in $\mathcal{K}$.
Now suppose that $\xi\in \mathcal{K}$, and then
  \begin{eqnarray*}
p(\xi)\,=\,\xi_{[0][0]}\, \int\Big(\xi_{[0][1]}\,S(\xi_{[1]})\Big) \,=\,
\xi_{[0]}\, \int\Big(\xi_{[1](1)}\,S(\xi_{[1](2)})\Big)\,=\,
\xi_{[0]}\, \epsilon(\xi_{[1]})\,=\,\xi\ .
\end{eqnarray*}
Finally we need to show the $H$-colinearity of $p$:
  \begin{eqnarray*}
p(\xi)_{[0]} \tens p(\xi)_{[1]}
  &=& p_0(\xi_{[0]})_{[0][0]} \tens
  p_0(\xi_{[0]})_{[0][1]}\, \int\Big(p_0(\xi_{[0]})_{[1]}\,S(\xi_{[1]})\Big) \cr
  &=& p_0(\xi_{[0]})_{[0]} \tens
  p_0(\xi_{[0]})_{[1]}\, \int\Big(p_0(\xi_{[0]})_{[2]}\,S(\xi_{[1]})\Big)\ ,\cr
  p(\xi_{[0]})\tens \xi_{[1]} &=&
  p_0(\xi_{[0][0]})_{[0]}\tens \xi_{[1]}\, 
\int\Big(p_0(\xi_{[0][0]})_{[1]}\,S(\xi_{[0][1]})\Big)\cr
  &=&
  p_0(\xi_{[0]})_{[0]}\tens \xi_{[2]}\, 
\int\Big(p_0(\xi_{[0]})_{[1]}\,S(\xi_{[1]})\Big)\ .
\end{eqnarray*}
Now Lemma~\ref{askfjauuu} gives the equality of these 
expressions.\quad$\square$

\begin{theorem}\label{bchsinbjic}
If $\pi:X\to H$ satisfies condition K and $H$ has a normalised left integral,
then, using the isomorphism \eqref{isom.2},
$\Omega^n B = (\mathcal{K}^{\wedge n}\tens X)^{{\rm co}H}$
(with the tensor coaction).
\end{theorem}
\noindent {\bf Proof:}\quad We shall prove this by induction on $n$, starting
at $\Omega^1 B = (\mathcal{K}\tens X)^{{\rm co}H}$, which has been done
in Corollary~\ref{askfjascv}.
Now assume the statement for $n$.

To show that $\Omega^{n+1} B \subset (\mathcal{K}^{\wedge n+1}\tens 
X)^{{\rm co}H}$,
we use $\Omega^{n+1} B\subset \Omega^{n} B\wedge \Omega^{1} B$, the formula
for the wedge product given in Section \ref{isuydfgc}, and 
Proposition~\ref{sjkhcv}.

To show $(\mathcal{K}^{\wedge n+1}\tens X)^{{\rm co}H}\subset \Omega^{n+1} B$,
take $(\xi\wedge\eta)\tens x\in (\mathcal{K}^{\wedge n+1}\tens X)^{{\rm co}H}$
(summation indices omitted for clarity), with $\xi\in L^1 X$
and $\eta\in (L^1 X)^{\wedge n}$. Then
\begin{eqnarray} \label{kjbasdcvb4}
(\xi\wedge\eta)\tens x \,=\,(p(\xi_{[0]})\tens S(\xi_{[1](1)}))\wedge 
(\xi_{[1](2)}\la
(\eta\tens x))\ ,
\end{eqnarray}
where we use square brackets for the right $X$-coaction on $L^1 X$ and
  $p:L^1 X\to L^1 X$ is the projection given in Lemma~\ref{askfjascvuuu}. In
  (\ref{kjbasdcvb4}) the first factor in the wedge product is in
  $\mathcal{K}\tens X$, and the second is in $\mathcal{K}^{\wedge n}\tens X$.
We use the colinearity of $p$ to rewrite
  (\ref{kjbasdcvb4}) as
  \begin{eqnarray} \label{kjbasdcvb5}
(\xi\wedge\eta)\tens x \,=\,(p(\xi_{[0][0]})\tens 
S(\xi_{[0][1]}))\wedge (\xi_{[1]}\la
(\eta\tens x))\ ,
\end{eqnarray}
and now it is evident that the first factor is in
$(\mathcal{K}\tens X)^{{\rm co}H}=\Omega^1 B$.
  It is not obvious that the second factor is $H$-invariant,
however we can integrate both sides of the equation over $H$, and this averages
the second factor to be $H$-invariant, without changing the left hand 
side. \quad $\square$

\section{Homogeneous spaces as fibrations}\label{sec.homog}
In this section we still consider Hopf algebras
$X$ and $H$ with bicovariant differential calculi. We assume  that there
exists a differentiable surjective Hopf algebra map $\pi:X\to H$. The 
differentiable right $H$-coaction
on $X$ is given by $\rho=(\id\tens\pi)\Delta:X\to X\tens H$.
\subsection{Checking the definition of fibration}

\begin{lemma}\label{kjhgc} Suppose that $H$ has a normalised left integral.
For any right $H$-comodule  and right $X$-module $V$ such that
the right action $\ra:V\tens X\to V$ is an $H$-comodule map
(with the tensor product coaction), the action
$\ra:V^{{\rm co}H}\tens_B X \to V$ is a bijective correspondence.
\end{lemma}
\noindent {\bf Proof:}\quad Call the coaction $\rho:V\to V\tens H$.
  Take a linear map $\psi:H\to X$ so that
$\pi\circ\psi=\id:H\to H$. Now define an inverse map
  $\tau:V\to V^{{\rm co}H}\tens_B X$ for $\ra:V^{{\rm co}H}\tens_B X \to V$ by
  \begin{eqnarray*}
\tau(a) \,=\,\Big((\id\tens\int)\rho\Big)\Big(a_{[0]}\ra 
S(\psi(a_{[1]})_{(1)}) \Big)
\tens_B \psi(a_{[1]})_{(2)}\ .
\end{eqnarray*}
The purpose of the operation $(\id\tens\int)\rho$ is to ensure that
the result lies in $V^{{\rm co}H}\tens_B X$ rather than just $V\tens_B X$.
The readers may verify that these two maps are inverse
by direct calculation. The reader familiar with the Hopf-Galois 
theory  may recognise this result as a consequence of 
\cite[Theorem~I]{Sch:pri}.
\quad$\square$

\begin{cor} \label{vcbhhasidj}
If $\pi:X\to H$ satisfies condition K, then, under the isomorphism 
\eqref{isom.2}, $\Omega^m B.X$
corresponds to $\mathcal{K}^{\wedge m}\tens X$.
\end{cor}
\noindent {\bf Proof:}\quad Put $V=\mathcal{K}^{\wedge m}\tens X$ in
Lemma \ref{kjhgc}, and use Theorem~\ref{bchsinbjic}.
\quad$\square$

\begin{cor} \label{uytre1}
For integer $n\ge 1$, $\Omega^m B\wedge\Omega^n X$
corresponds to $(\mathcal{K}^{\wedge m}\wedge(L^1 X)^{\wedge n})\tens X$,
and then
\begin{eqnarray*}
\Xi^0_m X &=& \Omega^m B.X \,=\, \mathcal{K}^{\wedge m} \tens X\ ,\cr
\Xi^n_m X &=& \frac{\Omega^m B\wedge\Omega^n X}
{\Omega^{m+1} B\wedge\Omega^{n-1} X} \cong
\frac{\mathcal{K}^{\wedge m}\wedge(L^1 X)^{\wedge n}}
{\mathcal{K}^{\wedge m+1}\wedge(L^1 X)^{\wedge n-1}} \tens X\ ,\quad n\ge 1\ .
\end{eqnarray*}
\end{cor}
\noindent {\bf Proof:}\quad Note that $\Omega^m B\wedge\Omega^n X =
\Omega^m B.X\wedge\Omega^n X$ and use Corollary \ref{vcbhhasidj}.
\quad$\square$

\begin{lemma}\label{uytre2}
  For integers $m\ge 0$ and $n\ge 1$, there is an isomorphism
\begin{eqnarray*}
\frac{\mathcal{K}^{\wedge m} \wedge (L^1 X)^{\wedge n}}
{\mathcal{K}^{\wedge m+1} \wedge (L^1 X)^{\wedge n-1}}\, \cong \,
\mathcal{K}^{\wedge m} \tens \frac{ (L^1 X)^{\wedge n}}
{\mathcal{K} \wedge (L^1 X)^{\wedge n-1}}
\end{eqnarray*}
induced by the map (for $\kappa_i\in \mathcal{K}$ and $\xi_i\in L^1 X$,
and  where $[.]$ denotes equivalence class)
\begin{eqnarray*}
[\kappa_1\tens\dots\tens\kappa_m\tens\xi_1\tens\dots\tens\xi_n]
\mapsto (\kappa_1\wedge\dots\wedge\kappa_m) \tens
[\xi_1\tens\dots\tens\xi_n]\ .
\end{eqnarray*}
\end{lemma}
\noindent {\bf Proof:}\quad First we must show that the map given in 
Lemma~\ref{uytre2} is well defined, and to do this we must use the 
braiding $\sigma$ in the definition of wedge product.
The left hand side space is defined to be the quotient of $\mathcal{K}^{\tens m}\tens
(L^1)^{\tens n}$  by a subspace spanned by elements of the form
$\kappa_1\tens\dots\tens\kappa_m\tens\xi_1\tens\dots\tens\xi_n$ where at least
one of the following is true:

(a)\quad $\xi_1\in\mathcal{K}$.

(b)\quad For some $1\le i\le m-1$, $\sigma(\kappa_i\tens\kappa_{i+1})=
\kappa_i\tens\kappa_{i+1}$.

(c)\quad $\sigma(\kappa_m\tens\xi_{1})=\kappa_m\tens\xi_{1}$.

(d)\quad For some $1\le i\le n-1$, $\sigma(\xi_i\tens\xi_{i+1})=
\xi_i\tens\xi_{i+1}$.

\noindent
We need to show that all these elements are mapped to zero. In cases
(a) and (d) we have $[\xi_1\tens\dots\tens\xi_n]=0$. In case
(b) we have $\kappa_1\wedge\dots\wedge\kappa_m=0$. In case (c),
\begin{eqnarray*}
\sigma(\kappa_m\tens\xi_{1})\,=\,\xi_{1[0]}\tens S(\xi_{1[1]})\la \kappa_m\,\in
\, L^1 X \tens \mathcal{K}\ ,
\end{eqnarray*}
since the left $X$ action restricts to $\mathcal{K}$. It then follows that
we are in case (a).

The inverse map is
\begin{eqnarray*}
[\kappa_1\tens\dots\tens\kappa_m]\tens [\xi_1\tens\dots\tens\xi_n] \mapsto
[\kappa_1\tens\dots\tens\kappa_m\tens\xi_1\tens\dots\tens\xi_n]\ ,
\end{eqnarray*}
and showing that this is well defined is rather easier than for the 
forward map.
\quad$\square$

\begin{theorem} \label{ikvyu}
Suppose that $X$ and $H$ are Hopf algebras with bicovariant 
differential structure,
and that $\pi:X\to H$ is a surjective differentiable Hopf
algebra map. Additionally suppose that

\smallskip
(1)\quad $\pi:X\to H$ satisfies condition K (see Definition \ref{fcond});

(2)\quad $H$ has a normalised left integral.

\smallskip
\noindent Then the inclusion $B=X^{{\rm co}H}\to X$ is a differentiable 
fibration (see Definition~\ref{jasdfvajkhv}), where
$\rho=(\id\tens\pi)\Delta_X:
X\to X\tens H$. Here $B$ has the differential structure given by 
Theorem~\ref{bchsinbjic}.
\end{theorem}
\noindent {\bf Proof:}\quad From Definition~\ref{jasdfvajkhv}
we need to show that the map $\Theta_m:\Omega^m B\tens_B \Xi_0^* X \to
\Xi_m^{*} X$ defined by
$\xi\tens [\eta]_0\mapsto [\xi\wedge\eta]_m$ is invertible for all $m\ge 0$.

We begin by using the fact that $\Omega^m B\tens_B \Xi_0^* X\cong
\Omega^m B\tens_B X \tens_X \Xi_0^* X$. Since we have $\Omega^m B=
(\mathcal{K}^{\wedge m}\tens X)^{{\rm co}H}$, Lemma~\ref{kjhgc} gives
$\Omega^m B\tens_B \Xi_0^* X\cong
(\mathcal{K}^{\wedge m}\tens X)\tens_X \Xi_0^* X$.
 From Section~\ref{isuydfgc} the right $X$-action on 
$\mathcal{K}^{\wedge m}\tens X$
is just multiplication on the second factor, so
\begin{eqnarray*}
\Omega^m B\tens_B \Xi_0^* X\,\cong\,
\mathcal{K}^{\wedge m}\tens X\tens_X \Xi_0^* X\,\cong\,
\mathcal{K}^{\wedge m}\tens  \Xi_0^* X\ ,
\end{eqnarray*}
and the result follows from
Corollary~\ref{uytre1} and Lemma~\ref{uytre2}.\quad$\square$

\subsection{Identifying the fibre of the fibration}
Assume the conditions of Theorem~\ref{ikvyu}

\begin{lemma} \label{kszcvbhj}
$L^1 X \wedge \mathcal{K}\subset \mathcal{K}\wedge L^1 X$.
\end{lemma}
\noindent {\bf Proof:}\quad From Corollary~\ref{vcbhhasidj},
$X.\extd B\subset X.\extd B.X=X.(\mathcal{K}\tens 
X)=(\mathcal{K}\tens X)= \extd B.X$.
Applying $\extd$ to this gives $\extd X\wedge \extd B \subset \extd 
B\wedge \extd X$.
 From this we conclude that  $\extd X\wedge \extd B.X \subset \extd 
B\wedge \extd X.X$,
so from condition K,
$\extd X\wedge \mathcal{K}\subset \extd B\wedge \extd X.X$. 
Multiplying again by $X$,
\begin{eqnarray*}
X.\extd X\wedge \mathcal{K}\subset X.\extd B\wedge \extd X.X \subset
\extd B\wedge X.\extd X.X\ ,
\end{eqnarray*}
so $\Omega^1 X \wedge \mathcal{K}\subset \extd B\wedge \Omega^1 X$.
 From Proposition~\ref{soidudvg}, $\extd B\subset \mathcal{K}.X$, so
$\Omega^1 X \wedge \mathcal{K}\subset \mathcal{K}\wedge \Omega^1 X$.
Then $L^1 X \wedge \mathcal{K}\subset \mathcal{K}\wedge \Omega^1 X$.
Note that $L^1 X \wedge \mathcal{K}$ is left invariant with respect to
the left $X$-coaction, and
consider $\xi\wedge\eta\in {}^{{\rm co}X} (\mathcal{K}\wedge \Omega^1 X)$.
Then, by invariance of $\xi\wedge\eta$ and $\xi$,
\begin{eqnarray*}
\xi\wedge\eta &=& S( (\xi\wedge\eta)_{[-1]}) (\xi\wedge\eta)_{[0]} \,=\,
S(\eta_{[-1]}) ( \xi\wedge\eta_{[0]}) \cr
&=& S(\eta_{[-1]})_{(3)} \xi \wedge S^{-1}(S(\eta_{[-1]})_{(2)}) \,
S(\eta_{[-1]})_{(1)} \,\eta_{[0]} \cr
&=& S(\eta_{[-1]}\sw 1)\, \xi\, S^{-1}(S(\eta_{[-1]}\sw 2))\wedge
  S(\eta_{[-1]}\sw 3) \,\eta_{[0]} \cr
  &=& S(\eta_{[-1]}\sw 1)\la \xi \wedge  S(\eta_{[-1]}\sw 2) \,\eta_{[0]}\cr
   &=& S(\eta_{[-1]})\la \xi \wedge  S(\eta_{[0][-1]}) \,\eta_{[0][0]}\ .
\end{eqnarray*}
As the left $X$-action restricts to $ \mathcal{K}$, this is in $ 
\mathcal{K}\tens L^1 X$.
\quad$\square$

\begin{propos}  The
  map $\pi_*^{\tens n}:(L^1 X)^{\tens n}\to (L^1 H)^{\tens n}$ induces an
invertible map
\begin{eqnarray*}
\tilde\pi:\frac{ (L^1 X)^{\wedge n}}
{\mathcal{K} \wedge (L^1 X)^{\wedge n-1}} \to (L^1 H)^{\wedge n}\ .
\end{eqnarray*}
\end{propos}
\noindent {\bf Proof:}\quad The domain of $\tilde\pi$ is the quotient of
$(L^1 X)^{\tens n}$  by a subspace spanned by elements of the form
$\xi_1\tens\dots\tens\xi_n$, where at least
one of the following is true:

(a)\quad $\xi_1\in\mathcal{K}$.

(b)\quad For some $1\le i\le n-1$, $\sigma(\xi_i\tens\xi_{i+1})=
\xi_i\tens\xi_{i+1}$.

\noindent The codomain of $\tilde\pi$ is the quotient of
$(L^1 H)^{\tens n}$  by a subspace spanned by elements of the form
$\eta_1\tens\dots\tens\eta_n$, where at least
one of the following is true:

(c)\quad For some $1\le i\le n-1$, $\sigma(\eta_i\tens\eta_{i+1})=
\eta_i\tens\eta_{i+1}$.

\noindent Case (a) maps to zero under $\pi_*^{\tens n}$ by definition of
$\mathcal{K}$. Case (b) maps to case (c) as $\pi$ is a Hopf algebra map.
Thus $\tilde\pi$ is well defined.

Given the hypothesis, it is automatic that $\tilde\pi$ is onto.

To show that $\tilde\pi$ is one-to-one, it is sufficient to show that the
subspace quotienting $(L^1 X)^{\tens n}$ in the domain contains all 
elements of the form
$\xi_1\tens\dots\tens\xi_n$, where at least
one of the following is true:

(d)\quad For some $1\le i\le n$, $\xi_i\in\mathcal{K}$.

\noindent This follows from repeated application of 
Lemma~\ref{kszcvbhj}. \quad$\square$

\section{Example: The non-commutative Hopf fibration with a 
bicovariant calculus}
In this section we return to the algebras $X$, $H$ and $B$ discussed 
in Section~\ref{kuasceewmnoon}, but now we consider a (minimal) 
bicovariant differential calculus on $\mathcal{A}(SL_q(2))$. In view 
of the results of Sections \ref{sec.bicov} and \ref{sec.homog} our 
task will be to construct a suitable calculus on $H$ so that the map 
$\pi:X\to H$ is differentiable and then to check that $\pi$ satisfies 
condition K in Definition~\ref{fcond}.

\subsection{A 4D bicovariant calculus on $\mathcal{A}(SL_q(2))$}
This differential calculus on $X=\mathcal{A}(SL_q(2))$ was introduced 
by Woronowicz in \cite{worondiff} and
is generated by four left invariant 1-forms
$\{\omega^1,\omega^2,\omega^+,\omega^-\}$.
  The differentials of the generators are given by
\begin{eqnarray} \label{4Drel1}
\extd \alpha &=& \frac{q-q^{-1}-q^{-2}}{q+1}\,\alpha\,\omega^1
- q^{-2}\,\beta\,\omega^++\frac{q^{-1}}{q+1}\,\alpha\,\omega^2\ ,\cr
\extd \beta &=& \frac{q}{q+1}\,\beta\,\omega^1
- q^{-2}\,\alpha\,\omega^- - \frac{q^{-2}}{q+1}\,\beta\,\omega^2\ ,\cr
\extd \gamma &=& \frac{q-q^{-1}-q^{-2}}{q+1}\,\gamma\,\omega^1
- q^{-2}\,\delta\,\omega^++\frac{q^{-1}}{q+1}\,\gamma\,\omega^2\ ,\cr
\extd \delta &=& \frac{q}{q+1}\,\delta\,\omega^1
- q^{-2}\,\gamma\,\omega^- - \frac{q^{-2}}{q+1}\,\delta\,\omega^2\ .
\end{eqnarray}
We have the commutation relations
\begin{eqnarray}\label{4Drel2}
\omega^2\,\alpha \,=\, 
q\,\alpha\,\omega^2-(q-q^{-1})\beta\,\omega^+
+q(q-q^{-1})^2\,\alpha\,\omega^1 &,&
\omega^2\,\beta \,=\, q^{-1}\beta\,\omega^2-(q-q^{-1})\alpha\,\omega^-\ ,\cr
\omega^-\,\alpha\,=\,\alpha\,\omega^--(q^2-1)\beta\,\omega^1 &,&
\omega^-\,\beta \,=\, \beta\,\omega^-\ ,\cr
\omega^+\,\alpha\,=\,\alpha\omega^+ &,&
\omega^+\,\beta\,=\,\beta\,\omega^+ -(q^2-1)\alpha\,\omega^1\ , \cr
\omega^1\,\alpha\,=\, q^{-1}\,\alpha\,\omega^1 &,&
\omega^1\,\beta\,=\,q\,\beta\omega^1\ .
\end{eqnarray}
and these relations with the replacements $\alpha\to \gamma$
and $\beta\to\delta$.

\subsection{The differentiability of $\pi:X\to H$}
We use the Hopf algebra map $\pi$ from Section~\ref{kuasceewmnoon},
and assume that the map $\pi:X\to H$ is differentiable, i.e.\ that it 
extends to
a map $\pi_*$ of differential graded algebras. Applying $\pi_*$ to 
the expression for
$\extd\beta$ in (\ref{4Drel1}) gives $z\,\pi_*(\omega^-)=0$, and 
since $z$ is invertible
we deduce $\pi_*(\omega^-)=0$. Likewise the expression for
$\extd\gamma$ in (\ref{4Drel1}) gives $\pi_*(\omega^+)=0$. Now the sixth
equation in (\ref{4Drel2}) gives $(q^2-1)\,z\,\pi_*(\omega^1)=0$, so if
$q\neq\pm 1$ we get $\pi_*(\omega^1)=0$. Now the equations
for $\extd\alpha$ and $\extd\delta$ in (\ref{4Drel1}) give
\begin{eqnarray*}
\extd z \,=\,\frac{q^{-1}}{q+1}\,z\,\pi_*(\omega^2)\ ,\quad
-z^{-1}.\extd z.z^{-1} \,=\,-\frac{q^{-2}}{q+1}\,z^{-1}\,\pi_*(\omega^2)\ .
\end{eqnarray*}
 From this we get
\begin{eqnarray}
\pi_*(\omega^2)\,=\,q(q+1)\,z^{-1}.\extd z\ ,\quad
\extd z.z^{-1}\,=\,q^{-1}\,z^{-1}.\extd z\ .
\end{eqnarray}
Just as in the case of the 3D calculus, we must have a noncommutative calculus
for the commutative algebra $H$. Note that $\mathcal{K}$, the left 
invariant forms
which are in the kernel of $\pi_*$, has basis $\{\omega^1,\omega^+,\omega^-\}$.

\subsection{Verifying condition K}\label{sub.check.F}

\begin{propos}
All of $\omega^-$, $\omega^+$ and $\omega^1$ are in $\extd B.X$.
\end{propos}
\noindent {\bf Proof:}\quad
Begin by calculating
\begin{eqnarray*}
\extd\beta_{(2)}.S^{-1}(\beta_{(1)}) &=& 
\extd\beta.\delta-q\,\extd\delta.\beta \cr
&=& \frac{q}{q+1}\,\beta\,\omega^1\,\delta
- q^{-2}\,\alpha\,\omega^-\,\delta - 
\frac{q^{-2}}{q+1}\,\beta\,\omega^2\,\delta \cr
&&-\, \frac{q^2}{q+1}\,\delta\,\omega^1\,\beta
+ q^{-1}\,\gamma\,\omega^-\,\beta + 
\frac{q^{-1}}{q+1}\,\delta\,\omega^2\,\beta \cr
&=&
- q^{-2}\,\alpha\delta\,\omega^- - \frac{q^{-2}}{q+1}\,\beta\,
(q^{-1}\delta\,\omega^2-(q-q^{-1})\gamma\,\omega^-) \cr
&&
+\, q^{-1}\,\gamma\beta\,\omega^- + \frac{q^{-1}}{q+1}\,\delta\,
(q^{-1}\beta\,\omega^2-(q-q^{-1})\alpha\,\omega^-) \cr
&=& - q^{-2}\,(\alpha\delta-q\,\gamma\beta)\,\omega^-
  - \frac{q^{-2}(q^2-1)}{q+1}\,(\delta\alpha-q^{-1}\beta\gamma)\,\omega^- 
  = -q^{-1}\,\omega^-\ ,
\end{eqnarray*}
and also
\begin{eqnarray*}
\extd\gamma_{(2)}.S^{-1}(\gamma_{(1)}) &=& 
\extd\gamma.\alpha-q^{-1}\,\extd\alpha.
\gamma \cr
&=& \frac{q-q^{-1}-q^{-2}}{q+1}\,\gamma\,\omega^1\,\alpha
- 
q^{-2}\,\delta\,\omega^+\,\alpha+\frac{q^{-1}}{q+1}\,\gamma\,\omega^2\,\alpha 
\cr
&&-\,q^{-1}\,\frac{q-q^{-1}-q^{-2}}{q+1}\,\alpha\,\omega^1\,\gamma
+ q^{-3}\,\beta\,\omega^+\,\gamma-\frac{q^{-2}}{q+1}\,\alpha\,\omega^2\,\gamma 
\cr
&=&
- q^{-2}\,\delta\alpha\,\omega^++\frac{q^{-1}}{q+1}\,\gamma
( q\,\alpha\,\omega^2-(q-q^{-1})\beta\,\omega^+
+q(q-q^{-1})^2\,\alpha\,\omega^1) \cr
&&
+\, q^{-3}\,\beta\gamma\,\omega^+-\frac{q^{-2}}{q+1}\,\alpha
(q\,\gamma\,\omega^2-(q-q^{-1})\delta\,\omega^+
+q(q-q^{-1})^2\,\gamma\,\omega^1) \cr
&=&
- q^{-2}\,\delta\alpha\,\omega^+-\frac{q^{-1}}{q+1}
  (q-q^{-1})\gamma\beta\,\omega^+ \cr
&&
+\, q^{-3}\,\beta\gamma\,\omega^++\frac{q^{-2}}{q+1}
(q-q^{-1})\alpha\delta\,\omega^+ \cr
&=&
- q^{-2}\,\omega^+
+\frac{q^{-3}(q^2-1)}{q+1}\,\omega^+\,=\,-q^{-3}\,\omega^+\ .
\end{eqnarray*}
Now we have
\begin{eqnarray*}
\extd(\alpha\beta)_{(2)}.S^{-1}((\alpha\beta)_{(1)}) &=&
\extd \alpha_{(2)}.S^{-1}(\alpha_{(1)})\,\epsilon(\beta) - q^{-1}\,
\alpha_{(2)}\,\omega^-\,S^{-1}(\alpha_{(1)}) \cr
&=& \gamma\,\omega^-\,\beta- q^{-1}\,\alpha\,\omega^-\,\delta \cr
&=& - q^{-1}\,(\alpha\delta-q\,\gamma\beta)\,\omega^-\,=\,
- q^{-1}\,\omega^-\ ,\cr
\extd(\gamma\beta)_{(2)}.S^{-1}((\gamma\beta)_{(1)}) &=&
\extd \gamma_{(2)}.S^{-1}(\gamma_{(1)})\,\epsilon(\beta) - q^{-1}\,
\gamma_{(2)}\,\omega^-\,S^{-1}(\gamma_{(1)}) \cr
&=& q^{-2}\,\alpha\,\omega^-\,\gamma-q^{-1}\,\gamma\,\omega^-\,\alpha \cr
&=& q^{-2}\,\alpha\,(\gamma\,\omega^--(q^2-1)\delta\,\omega^1)
-q^{-1}\,\gamma\,(\alpha\,\omega^--(q^2-1)\beta\,\omega^1) \cr
&=& q^{-1}\, (q^2-1)\,\gamma\beta\,\omega^1
-q^{-2}\, (q^2-1)\,\alpha\delta\,\omega^1\cr
&=&  (q^{-2}-1)\,(\alpha\delta-q\,\beta\gamma)\,\omega^1
\,=\,  (q^{-2}-1)\,\omega^1\ ,
\cr
\extd(\delta\gamma)_{(2)}.S^{-1}((\delta\gamma)_{(1)}) &=&
\extd \delta_{(2)}.S^{-1}(\delta_{(1)})\,\epsilon(\gamma) - q^{-3}\,
\delta_{(2)}\,\omega^+\,S^{-1}(\delta_{(1)}) \cr
&=& q^{-4}\,\beta\,\omega^+\,\gamma-q^{-3}\,\delta\,\omega^+\,\alpha \cr
&=& -q^{-3}\,(\delta\alpha-q^{-1}\,\beta\gamma)\,\omega^+
\,=\, -q^{-3}\,\omega^+\ .
\end{eqnarray*}
This proves the claim since, for all $b\in B$, $\Delta b\in X\otimes 
B$. $\quad\square$

\bigskip \noindent {\bf \large Acknowledgements}\quad The authors 
would like to thank 
F.W.\ Clarke and M.D.\ Crossley (Swansea) for 
their help. EJB would like to thank Rutgers University,
where part of 
this work was done while he was on leave from Swansea, and the New 
York Public Libraries for their help. The research of TB is also supported by
the EPSRC grant GR/S01078/01.

\end{document}